\documentclass[12pt,fleqn]{amsart}
\usepackage{geometry}
\usepackage[T1]{fontenc}
\usepackage[
  style=alphabetic,
  maxalphanames=99,
]{biblatex}
\renewbibmacro*{volume+number+eid}{%
  \printfield{volume}%
  \printfield{number}%
  \setunit{\addcomma\space}%
  \printfield{eid}}
\DeclareFieldFormat[article]{number}{\mkbibparens{#1}}

\usepackage{mathtools}
\usepackage{thmtools}
\usepackage{amsmath}
\usepackage{amssymb}
\usepackage{mathrsfs}
\usepackage{array}

\usepackage{hyperref}
\hypersetup{
    urlcolor=blue,
}

\usepackage{enumerate}
\usepackage[shortlabels]{enumitem}
\usepackage{multicol}

\usepackage{setspace}

\usepackage{graphicx}

\usepackage{changepage}
\usepackage{wrapfig}

\usepackage{pifont}

\usepackage{cellspace}
\usepackage{listings}

\usepackage[normalem]{ulem}

\usepackage{outlines}

\usepackage{tikz}
\usepackage{tikz-cd}

\usepackage{subfig}

\usepackage[only,llbracket,rrbracket]{stmaryrd}
\DeclarePairedDelimiter\bbrackets{\llbracket}{\rrbracket}

\usepackage{needspace}

\usepackage{amsthm}
\counterwithin{figure}{section}

\usepackage[capitalise]{cleveref}

\theoremstyle{definition}
\newtheorem{theorem}{Theorem}[section]
\newtheorem{definition}[theorem]{Definition}
\newtheorem{lemma}[theorem]{Lemma}
\newtheorem{corollary}[theorem]{Corollary}

\newtheorem{construction}[theorem]{Construction}

\theoremstyle{remark}
\newtheorem{remark}[theorem]{Remark}

\DeclareMathOperator{\Hom}{Hom}

\DeclareMathOperator{\Clp}{Clp}
\DeclareMathOperator{\Uf}{Uf}

\newcommand{\bb}[1]{\mathbf{#1}}
\newcommand{\frk}[1]{\mathfrak{#1}}
\newcommand{\mca}[1]{\mathcal{#1}}
\newcommand{\msf}[1]{\mathsf{#1}}

\newcommand{\setbuilder}[2]{\left\{#1 : #2\right\}}

\newcommand{\set}[1]{\left\{#1\right\}}

\newcommand{\tand}{\text{ and }}

\newcommand{\la}{\langle}
\newcommand{\ra}{\rangle}

\DeclarePairedDelimiter{\abs}{\lvert}{\rvert}

\title{Local tabularity in MS4 with Casari's axiom}

\author{Chase Meadors}
\address{Department of Mathematics\\
University of Colorado Boulder\\
Boulder, CO, 80309\\
USA}
\email{chase.meadors@colorado.edu}


\DeclareMathOperator{\qmax}{qmax}

\newcommand{\mpsfl}{\mathsf{M^{+}S4}}
\newcommand{\mpsfv}{\mathbf{M^{+}S4}}
\newcommand{\mpgrzl}{\mathsf{M^{+}Grz}}

\newcommand{\approxexists}[1]{#1^{\approx}_{\exists}}

\begin{document}

\begin{abstract}
  We study local tabularity (local finiteness) in some extensions of $\msf{MS4}$ (monadic $\msf{S4}$).
  Our main result is a semantic characterization of local finiteness in varieties of $\mpsfl$-algebras, where $\mpsfl$ denotes the extension of $\msf{MS4}$ by the Casari axiom.
  We improve this to a syntactic criterion via the \textit{reducible path property} identified in \cite{Shap16}, and note that the product logic $\msf{S4}[n] \times \msf{S5}$ is an extension of $\mpsfl$, obtaining a criterion for extensions of $\msf{S4}[n] \times \msf{S5}$ as an application.
  Next, we give a characterization of local finiteness in varieties of $\msf{MS4B}[2]$-algebras, where $\msf{MS4B}$ denotes the extension of $\msf{MS4}$ by the Barcan axiom.
  We demonstrate that our methods cannot be extended beyond depth 2, as we give a translation of the fusion $\msf{S5}_2$ into $\msf{MS4B}[3]$ for $n \geq 3$ that preserves and reflects local finiteness, suggesting that a characterization there remains difficult.
  Finally, we also establish the finite model property for some of these logics which are not known to be locally tabular.
\end{abstract}

\maketitle

\tableofcontents

\section{Introduction}

It is a classic result in modal logic, known as the Segerberg--Maksimova theorem, that a normal extension of $\msf{S4}$ is locally tabular iff it is of finite depth (\cite{Seg71}, \cite{Mak75}).
Since a logic is locally tabular iff its corresponding variety is locally finite, this provides a characterization of locally finite varieties of $\msf{S4}$-algebras.

An important extension of $\msf{S4}$ is the bimodal logic $\msf{MS4}$---monadic $\msf{S4}$---which axiomatizes the one-variable fragment of predicate $\msf{S4}$ (\cite{FS77}).
\cite{BM23} addresses the question of how and to what extent the Segerberg--Maksimova theorem can be generalized to the monadic setting.
The main result of \cite{BM23} is negative, demonstrating a translation of the fusion $\msf{S5}_2$ into a rather strong extension of $\msf{MS4}$; nevertheless, some families of extensions of $\msf{MS4}$ were identified for which a generalized version of Segerberg--Maksimova holds.
This paper obtains more positive results in this direction.
One main result is a semantic and syntactic criterion for local tabularity, directly generalizing Segerberg--Maksimova, in $\mpsfl$, the extension of $\msf{MS4}$ obtained by asserting the G\"odel translation of the Casari axiom (\cref{thm:mpsfl-syntactic-criterion}).
This result contributes to a body of evidence that $\mpsfl$ is a more natural and well-behaved alternative to $\msf{MS4}$; for example, $\mpsfl$ plays a central role in obtaining a faithful provability interpretation of monadic intuitionistic predicate logic, as shown in \cite{BBI23}.
Our syntactic characterization is via the \textit{reducible path property} identified in \cite{Shap16}; we also provide an application, demonstrating that this yields a criterion for extensions of $\msf{S4}[n] \times \msf{S5}$ since, as the authors of \cite{Shap16} observed, this logic is an extension of $\mpsfl$.

We denote the extension of $\msf{MS4}$ by the \textit{Barcan axiom} by $\msf{MS4B}$.
The other main result of the paper is a syntactic criterion for local tabularity for extensions of $\msf{MS4B}[2]$ in \cref{thm:ms4b-syntactic-criterion}, extending a result of Shapirovsky and Sliusarev for extensions of the product logic $\msf{S4.1}[2] \times \msf{S5}$ (\cite[Thm.~6.10]{Shap24}).
However, we demonstrate that such a characterization when $n \geq 3$ seems out of reach, by showing that our translation of $\msf{S5}_2$ from \cite{BM23} can be modified to land in $\msf{MS4B}[3]$.

\cref{sec:prelims} gives the relevant preliminary notions concerning algebraic and relational semantics for $\msf{MS4}$.
In \cref{sec:fmp}, we establish the finite model property for some extensions of $\msf{MS4}$ that are not locally tabular.
\cref{sec:casari}, we develop the theory of $\mpsfl$-algebras and their dual frames, culminating in the characterization of local tabularity (\cref{thm:mpsfl-syntactic-criterion}).
We then demonstrate as a consequence that this criterion also characterizes local tabularity above $\msf{S4}[n] \times \msf{S5}$.
Likewise, \cref{sec:products} develops the corresponding theory for $\msf{MS4B}$, culminating in the characterization of local tabularity above $\msf{MS4B}[2]$ (\cref{thm:ms4b-syntactic-criterion}), and then covers the aforementioned translation of $\msf{S5}_2$ into $\msf{S4}[3] \times \msf{S5}$.
Finally, \cref{sec:conclusion} discusses possible future directions of research.

\section{Preliminaries}
\label{sec:prelims}

\subsection{\texorpdfstring{$\msf{MS4}$}{MS4}, algebraic semantics, and duality}

Let $\mathcal L$ be a propositional modal language with two modalities $\lozenge$ and $\exists$. As usual, we write $\square = \neg \lozenge \neg$ and $\forall = \neg \exists \neg$.

\begin{definition}[\cite{FS77}]
    The bimodal logic $\mathsf{MS4}$ is the smallest normal modal logic in $\mathcal L$ containing 
\begin{itemize}
    \item The $\mathsf{S4}$ axioms for $\lozenge$ (i.e.~the $\msf{K}$ axiom along with $p \to \lozenge p$ and $\lozenge \lozenge p \to \lozenge p$),
    \item The $\mathsf{S5}$ axioms for $\exists$ (e.g.~the $\msf{S4}$ axioms along with $\exists \forall p \to \forall p$),
    \item The left commutativity axiom $\exists \lozenge p \to \lozenge \exists p$.
\end{itemize}
\end{definition}

\begin{remark}
    In the terminology of \cite{GKWZ03}, $\msf{MS4} = [\msf{S4}, \msf{S5}]^\text{EX}$; that is, $\msf{MS4}$ is the \textit{expanding relativized product} of $\msf{S4}$ and $\msf{S5}$.
\end{remark}

Algebraic semantics for $\msf{MS4}$ is given by the following Boolean algebras with operators, first considered by Fischer-Servi \cite{FS77} under the name of \textit{bimodal algebras}.

\begin{definition}
\label{def:ms4}
An $\msf{MS4}$-algebra is a tuple $\frk{A} = (B, \lozenge, \exists)$ such that
\begin{itemize}
    \item $(B, \lozenge)$ is an $\msf{S4}$-algebra, i.e.~$B$ is a Boolean algebra and $\lozenge$ is a unary function on $B$ satisfying the identities of a closure operator:
    \[
    \lozenge 0 = 0 \qquad \lozenge (a \vee b) = \lozenge a \vee \lozenge b \qquad a \leq \lozenge a \qquad \lozenge \lozenge a \leq \lozenge a.
    \]
    \item $(B, \exists)$ is an $\msf{S5}$-algebra, that is an $\msf{S4}$-algebra that in addition satisfies $\exists \forall a \leq \forall a$ (where we write $-$ for Boolean negation).
    \item $\frk{A}$ satisfies $\exists \lozenge a \leq \lozenge \exists a$ (the \textit{left commutativity} axiom or simply the $\msf{MS4}$-axiom).
\end{itemize}
\end{definition}

It is clear that the class of $\msf{MS4}$-algebras is equationally definable, and hence forms a variety. We denote it and the corresponding category by $\bb{MS4}$.

We will make use of J\'{o}nsson--Tarski duality \cite{JT51} to work with the dual spaces of $\msf{MS4}$-algebras. We recall some relevant notions:

\begin{definition} 
    A topological space $X$ is a \textit{Stone space} if it is compact, Hausdorff, and zero-dimensional---that is, $X$ has a basis of {\em clopen sets} (sets that are both closed and open).
\end{definition}

\begin{definition} \label{def: cont rel}
    Let $X$ be a Stone space.
    We say a relation $R \subseteq X^2$ is \textit{continuous} if 
    \begin{enumerate}
        \item $R(x) := \setbuilder{y \in X}{x R y}$ is closed for each $x \in X$ ($R$ is \textit{point-closed}), and
        \item $R^{-1}(U) := \setbuilder{y \in X}{yRx \text{ for some } x \in U}$ is clopen whenever $U \subseteq X$ is clopen.
    \end{enumerate}
\end{definition}

Specializing multimodal descriptive general frames (see, e.g., \cite[Sec.~5.5]{BRV01}) to $\msf{MS4}$ yields: 

\begin{definition}
A {\em descriptive $\msf{MS4}$-frame} is a tuple $\frk{F} = (X, R, E)$ such that $X$ is a Stone space, $R$ is a quasi-order (reflexive and transitive), and $E$ is an equivalence relation on $X$ such that
\begin{itemize}
    \label{def:ms4-frame}
    \item both $R$ and $E$ are continuous relations, 
    \item $RE \subseteq ER$---i.e., $\forall x,y,y' \in X \; (x E y \tand y R y') \to \exists x' \in X \; (x R x' \tand x' E y')$.

    \begin{center}
    \begin{tikzpicture}[
        scale=1.5,
    	dot/.style={circle,fill=black,minimum size=6pt,inner sep=0}
    ]
    \node[dot,label=left:$x$] (bl) at (0, 0) {};
    \node [dot,label=left:$x'$] (tl) at (0, 1) {};
    \node [dot,label=right:$y'$] (tr) at (1, 1) {};
    \node [dot,label=right:$y$] (br) at (1, 0) {};
    \draw[latex-latex] (bl) -- (br);
    \draw[latex-latex,dashed] (tl) -- (tr);
    \draw[-latex,dashed] (bl) -- (tl);
    \draw[-latex] (br) -- (tr);
    \node at (0.5, -0.2) {$E$};
    \node at (-0.2, 0.5) {$R$};
    \node at (0.5, 1.2) {$E$};
    \node at (1.2, 0.5) {$R$};
    \end{tikzpicture}
    \end{center}
\end{itemize}
\end{definition}

We will refer to descriptive frames simply as \textit{frames} because they are the only kind of frames we will deal with in this paper (as opposed to Kripke frames or the broader class of general frames studied in modal logic \cite{CZ97}, \cite{BRV01}).

The next definition is well known (see, e.g., \cite[Sec.~3.3]{BRV01}).

\begin{definition} 
    Let $X, X'$ be sets, $S$ a binary relation on $X$, and $S'$ a binary relation on $X'$.
    A function $f : X \to X'$ is a \textit{p-morphism} (or \textit{bounded morphism})
    with respect to $(S, S')$ if 
    \begin{itemize}
        \item $x S y$ implies $f(x) S' f(y)$;
        \item $f(x) S' y'$ implies $x S y$ for some $y \in X$ with $f(y) = y'$.
    \end{itemize}
\end{definition}

\begin{definition}
    \label{def:ms4-p-morphism}
    An {\em $\msf{MS4}$-morphism} is a continuous map $f:(X, R, E) \to (X', R', E')$ between $\msf{MS4}$-frames such that $f$ is a p-morphism with respect to both $(R, R')$ and $(E, E')$.
\end{definition}

Specializing J\'{o}nsson--Tarski duality 
to $\msf{MS4}$-algebras yields:

\begin{theorem}
    The category of $\msf{MS4}$-algebras with homomorphisms and the category of descriptive $\msf{MS4}$-frames with $\msf{MS4}$-morphisms are dually equivalent.
\end{theorem}

\begin{remark}
\label{rem:jt-duality}
The above dual equivalence is implemented 
as follows:
with each $\msf{MS4}$-frame $\mathfrak{F} = (X, R, E)$ we associate the $\msf{MS4}$-algebra $\mathfrak{F}^* = (\Clp X, \lozenge_R, \lozenge_E)$ of clopen subsets of $X$, where $\lozenge_R$ and $\lozenge_E$ are the \textit{dual operators} of $R$ and $E$, defined by
\[
\lozenge_R U = R^{-1}(U) \ \mbox{ and } \ \lozenge_E U = E(U).
\]
In the other direction, with each $\msf{MS4}$-algebra $\frk{A} = (B, \lozenge, \exists)$ we associate the descriptive frame $\frk{A}_* = (\Uf B, R_\lozenge, R_\exists)$ of ultrafilters of $B$, where $R_\lozenge$ and $R_\exists$ are the \textit{dual relations} of $\lozenge$ and $\exists$, defined by 
\[
x R_\lozenge y \mbox{ iff } \forall a, a \in y \to \lozenge a \in x \ \mbox{ and } \ x R_\exists y \mbox{ iff } \forall a, \exists a \in x \leftrightarrow \exists a \in y,
\]
and $\Uf B$ is equipped with the \emph{Stone topology} generated by the clopen basis $\setbuilder{\frk{s}(a)}{a \in B}$, where $\frk{s}(a) = \setbuilder{x}{a \in x}$.
The unit isomorphisms of the duality are given by
\begin{itemize}
    \item $\frk{A} \to (\frk{A}^*)_*$; \ $a \mapsto \setbuilder{x}{a \in x}$,
    \item $\frk{F} \to (\frk{F}^*)_*$; \ $x \mapsto \setbuilder{a}{x \in a}$.
\end{itemize}
The dual of a morphism $f$ (in both categories) is the inverse image map $f^{-1}[\cdot]$, and   
the natural bijection $\Hom(\frk{A}, \frk{F}^*) \cong \Hom(\frk{F}, \frk{A}_*)$ is given by associating to $f : \frk{A} \to \frk{F}^*$ the morphism $\bar{f} : \frk{F} \to \frk{A}_*$ defined by $x \mapsto \setbuilder{a \in \frk{A}}{x \in f(a)}$, and to $g : \frk{F} \to \frk{A}_*$ the morphism $\bar{g} : \frk{A} \to \frk{F}^*$ defined by $a \mapsto \setbuilder{x \in \frk{F}}{a \in f(x)}$.
\end{remark}

In the sequel, we will make heavy use of this duality, often identifying an $\msf{MS4}$-algebra $\frk{A}$ with the algebra of clopen subsets of its dual frame $\frk{A}_*$, interpreting the Boolean connectives and operators as described in \cref{rem:jt-duality}.
For a formula $\varphi$ and an $\msf{MS4}$-algebra $\frk{A}$, we write $\frk{A} \models \varphi$ to mean $\frk{A} \models \varphi \approx 1$ in the usual sense of universal algebra (conflating a formula with its corresponding identity). 
It is well-known (see, e.g., \cite[Thm.~6.3]{Ven07}) that the lattice of subvarieties of $\bb{MS4}$ is dually isomorphic to the lattice of normal extensions of $\msf{MS4}$, hence $\msf{MS4}$-algebras provide adequate algebraic semantics for $\msf{MS4}$.
On the other hand, descriptive $\msf{MS4}$-frames may be understood as providing adequate relational/Kripke semantics for $\msf{MS4}$, and accordingly the notion $\frk{F} \models \varphi$ for $\frk{F} = (X, R, E)$ an $\msf{MS4}$-frame can be defined internally.
Indeed, for a valuation $\nu : \msf{Var} \to \Clp X$ assigning propositional variables to clopen sets and a point $x \in X$, define $(\frk{F}, x) \models_\nu \varphi$ (or simply $x \models \varphi$, when $\frk{F}$ and $\nu$ are clear from context) by induction:
\[
    \begin{array}{lcll}
        x \models p & \Leftrightarrow & x \in \nu(p) & (p \in \msf{Var}) \\
        x \models \varphi \cdot \psi & \Leftrightarrow & (x \models \varphi) \cdot (x \models \psi) & (\cdot \in \set{\wedge, \vee}) \\
        x \models \neg \varphi & \Leftrightarrow & \neg(x \models \varphi) & \\
        x \models \lozenge \varphi & \Leftrightarrow & \exists y, xRy \wedge y \models \varphi & \\
        x \models \exists \varphi & \Leftrightarrow & \exists y, xEy \wedge y \models \varphi
    \end{array}
\]
This definition extends $\nu$ to a map $\bbrackets{\cdot}_\nu : \msf{Fml} \to \Clp X$ by $\bbrackets{\varphi}_\nu = \setbuilder{x \in X}{x \models_\nu \varphi}$ (again we will drop the subscript $\nu$ when clear from context).
Now define $\frk{F} \models \varphi$ to mean $\bbrackets{\varphi}_\nu = X$ for all valuations $\nu$.
It is clear that $\frk{F} \models \varphi$ iff $\frk{F}^* \models \varphi$, but we will have occasion to use relational/frame-theoretic reasoning as opposed to algebraic reasoning (e.g., \cref{thm:cas-frame-condition} and \cref{lem:minimal-variety-axiomitization}).

\subsection{Depth of \texorpdfstring{$\msf{MS4}$}{MS4}-frames and related notions}
\label{sec:prelims-depth}

Describing the dual spaces of $\mpsfv$-algebras (the content of \cref{sec:casari-dual-spaces}) requires several concepts and some associated notation that we discuss here.

Let $\frk{A}$ be an $\msf{S4}$-algebra and $\frk{F}$ its dual frame.
Let $R^+ = R - R^{-1}$, i.e. $R^+$ is the relation of being \textit{properly} $R$-related.
By the \textit{depth} of a point $x \in \frk{F}$, we mean the supremum of the lengths of $R^+$-chain emanating from $x$, i.e. the largest $n$ such that there exists $x = x_1 \, R^+ \, x_2 \, R^+ \dots R^+ \, x_n$, in which case we write $d(x) = n$, or if there are arbitrarily long chains emanating from $x$, $d(x) = \omega$.
The depth of an algebra $\frk{A}$ or its dual frame $\frk{F}$ is $d(\frk{A}) = d(\frk{F}) = \sup \setbuilder{d(x)}{x \in \frk{F}}$.
For a variety $\bb{V} \subseteq \bb{S4}$, we write $d(\bb{V}) \leq n$ if for all $\frk{A} \in \bb{V}$, $d(\frk{A}) \leq n$, and $d(\bb{V}) = n$ if moreover there exists some $\frk{A} \in \bb{V}$ with $d(\frk{A}) = n$.
If $\bb{V}$ contains algebras of arbitrary depth, we say $d(\bb{V}) = \omega$.
We will also apply these notions to (varieties of) $\msf{MS4}$-algebras and their dual frames, understanding the depth of an $\msf{MS4}$-algebra to be the depth of its $\msf{S4}$-reduct (that is, the $R$-depth of its dual frame).

Consider the family of formulas $P_n$ defined by
\[
P_1 = \lozenge \square q_1 \to \square q_1 \qquad P_n = \lozenge (\square q_n \wedge \neg P_{n-1}) \to \square q_n.
\]
It is well known (see, e.g., \cite[Thm.~3.44]{CZ97}) that a variety $\bb{V} \subseteq \bb{S4}$ has $d(\bb{V}) \leq n$ iff $\bb{V} \models P_n$.
Following \cite{Shap16}, we write $\bb{V}[n]$ to refer to the subvariety $\bb{V} + P_n$.

Closely related to this are the notions of quasi-maximal points and the \textit{layers} of an $\msf{MS4}$-frame:

\begin{definition}
    Let $\frk{F} = (X, R, E)$ be an $\msf{MS4}$-frame.
    For $A \subseteq X$, we write 
    \[
        \qmax A  = \setbuilder{x \in A}{x R y \text{ for } y \in A \text{ implies } y R x}
    \]
    and define the \textit{$n$-th layer} of $\frk{F}$ inductively as
    \[
        D_0(\frk{F}) = \varnothing, \qquad D_{n+1}(\frk{F}) = \qmax \left( X - \bigcup_{i=0}^n D_i(\frk{F}) \right)
    \]
    (in particular, $D_1(\frk{F}) = \qmax X$).
    We say that a point $x \in X$ is \textit{of depth $n$} if $x \in D_n(\frk{F})$ and \textit{of finite depth} if $x \in D_n(\frk{F})$ for some $n < \omega$.
\end{definition}

For an arbitrary quasi-order, some or all $D_i$ may be empty (e.g., when $R$ contains infinite ascending chains).
In the case of descriptive $\msf{MS4}$-frames however, the following well-known results of Esakia apply:

\begin{theorem}[{\cite[Thm.~3.2.1-3]{Esa19}}]
    \label{thm:sees-maximal}
    Let $\frk{F} = (X, R, E)$ be an $\msf{MS4}$-frame and $U \subseteq X$.
    \begin{enumerate}
        [label=\normalfont(\arabic*), ref = \thetheorem(\arabic*)]
        \item \label[theorem]{thm:sees-maximal:1} If $U$ is closed, then for every $x \in U$ there exists $y \in \qmax U$ with $xRy$.
        \item \label[theorem]{thm:sees-maximal:2} $\qmax X$ is non-empty and closed.
    \end{enumerate}
\end{theorem}

In particular this implies that $D_1(\frk{F})$ is non-empty and closed for any descriptive $\msf{MS4}$-frame $\frk{F}$.
Clearly if $\frk{F}$ is of finite depth $n$, then $D_1(\frk{F}), \dots, D_n(\frk{F})$ are precisely the non-empty layers of $\frk{F}$.

Every $\msf{MS4}$-algebra (resp.~$\msf{MS4}$-frame) has a naturally associated $\msf{S4}$-algebra (resp.~$\msf{S4}$-frame).

\begin{definition}
    \label{def:b0}
    Let $\frk{A} = (B, \lozenge, \exists)$ be an $\msf{MS4}$-algebra and $\frk{F} = (X, R, E)$ its dual frame.
    \begin{itemize}
        \item Let $\frk{A}^0 = (B^0, \lozenge^0)$, where $B^0 = \exists B$ is the set of $\exists$-fixpoints of $B$, and $\lozenge^0$ the restriction of $\lozenge$ to $B^0$.
        \item Let $\frk{F}^0 = (X/E, \overline{R})$, where $\alpha \overline{R} \beta$ iff $\exists x \in \alpha, y \in \beta, x R y$.
    \end{itemize}
\end{definition}

It is a direct consequence of the $\msf{MS4}$ axiom that $\frk{A}^0$ is an $\msf{S4}$-subalgebra of $(B, \lozenge)$ and, dually, $\frk{F}^0$ is a continuous p-morphic image of $(X, R)$; indeed, $\frk{F}^0$ is precisely the dual frame of the subalgebra $\frk{A}^0$.

\begin{remark}
    \label{rem:ms4-s4-representation}
    Similar to monadic Heyting algebras (see, e.g., \cite[Sec.~3]{Bez98}), 
    we have that $\msf{MS4}$-algebras can be represented as pairs of $\msf{S4}$-algebras $(\frk{A}, \frk{A}_0)$ such that $\frk{A}_0$ is an $\msf{S4}$-subalgebra of $\frk{A}$ and the inclusion $\frk{A}_0 \hookrightarrow \frk{A}$ has a left adjoint $(\exists)$.
\end{remark}

$\frk{A}^0$ and $\frk{F}^0$ are closely related to the following auxiliary notions:

\begin{definition}\
    \begin{enumerate}
        \item For an $\msf{MS4}$-algebra $\mathfrak A=(B, \lozenge, \exists)$, define $\blacklozenge = \lozenge \exists$.
        \item For an $\msf{MS4}$-frame $\mathfrak F=(X, R, E)$, define $Q = ER$.
    \end{enumerate}
\end{definition}

Clearly $\blacklozenge$ is an $\bb{S4}$-possibility operator on $B$ and $Q$ is a quasi-order on $X$; indeed, $Q$ is the dual relation of $\blacklozenge$ (see \cite[Lem.~3.2]{BM23}).

Since the universe of $\frk{A}^0$ consists of the $\exists$-fixpoints of $\frk{A}$, we have that $\frk{A}^0 \models \varphi$ iff $\frk{A} \models \varphi^0$, where $\varphi^0$ denotes the formula obtained by replacing all occurrences of $\lozenge$ in $\varphi$ with $\blacklozenge$.
As a consequence, the $\blacklozenge$-depth of $\frk{A}$ and the $\lozenge^0$-depth of $\frk{A}^0$ coincide and, dually, the $Q$-depth of $\frk{F}$ and the $\overline{R}$-depth of $\frk{F}^0$ coincide.
For this reason, we introduce the following notation:

\begin{definition}
    Let $\frk{A} = (B, \lozenge, \exists)$ be an $\msf{MS4}$-algebra and $\frk{F} = (X, R, E)$ its dual frame.
    \begin{itemize}
        \item Define $d^0(\frk{A}) = d(\frk{A}^0)$ and $d^0(\frk{F}) = d(\frk{F}^0)$.
        \item Let $P^0_n$ be the formula $P_n$ with all instances of $\lozenge$ replaced by $\blacklozenge$.
    \end{itemize}    
\end{definition}

It is immediate that $\bb{V} \models P^0_n$ iff $d^0(\bb{V}) \leq n$.
Moreover, these notions are precisely the ones necessary to characterize the dual frames of simple and subdirectly irreducible (hereafter abbreviated s.i.) $\msf{MS4}$-algebras.

\begin{definition}
    Let $\frk{F} = (X, R, E)$ be an $\msf{MS4}$-frame.
    \begin{itemize}
        \item $x \in X$ is a \textit{$Q$-root} of $\frk{F}$ if $Q(x) = X$.
        \item $T(\frk{F})$ denotes the set of $Q$-roots of $\frk{F}$.
        \item $\frk{F}$ is \textit{$Q$-rooted} if $T(\frk{F})$ is non-empty, and \textit{strongly $Q$-rooted} if moreover $T(\frk{F})$ is open.
    \end{itemize}
\end{definition}

\begin{theorem}[{\cite[Thm.~3.5]{BM23}}]
    Let $\frk{A}$ be a (nontrivial) $\msf{MS4}$-algebra and $\frk{F} = (X, R, E)$ its dual $\msf{MS4}$-frame.
    \begin{enumerate}
        [label=\normalfont(\arabic*), ref = \thetheorem(\arabic*)]
        \item \label[theorem]{thm:si-simple-ms4:1} $\frk{A}$ is s.i.~iff $\frk{F}$ is strongly $Q$-rooted.
        \item \label[theorem]{thm:si-simple-ms4:2} $\frk{A}$ is simple iff $T(\frk{F}) = X$.
    \end{enumerate}
\end{theorem}

We denote by $\bb{V}_{\text{SI}}$ the class of subdirectly irreducible $\bb{V}$-algebras.
Recall that a variety is \textit{semisimple} if all its subdirectly irreducible members are simple.
We define $\bb{MS4_S} = \bb{MS4} + P^0_1$ and have

\begin{theorem}[{\cite[Thm.~3.8]{BM23}}]
    \label{thm:ms4s-semisimple}
    $\bb{MS4_S}$ is the largest semisimple subvariety of $\bb{MS4}$.
\end{theorem}

\section{The fmp for some subvarieties of \texorpdfstring{$\bb{MS4}$}{MS4} of bounded depth}
\label{sec:fmp}

In \cite{BM23}, the finite model property was established for $\bb{MS4_S}$, the largest semisimple subvariety of $\bb{MS4}$.
The varieties $\bb{MS4}[n]$ and even $\bb{MS4_S}[n]$ are not locally finite, as they both contain the variety $\bb{MS4}[1] = \bb{S5}^2$ of two-dimensional diagonal-free cylindric algebras, which is not locally finite by \cite[Thm.~2.1.11]{HMT85} (see also \cite[Ex.~5.10]{BM23}).
Furthermore, the question of characterizing local finiteness in these varieties is expected to be difficult to resolve, in light of \cite[Sec.~6]{BM23}.
As such, it was left open whether these natural subvarieties of $\bb{MS4}$ have the finite model property.
This section resolves the question in the positive for both.

\begin{definition}
  Let $\frk{A} = (B, \lozenge, \exists)$ be an $\msf{MS4}$-algebra, and let $B' \subseteq B$ be a finite $\msf{S4}$-subalgebra of $(B, \lozenge)$.
  Define $\exists' : B' \to B'$ by
  \[
    \exists' a = \bigwedge \setbuilder{x \in B^0 \cap B'}{a \leq x}
  \]
  And define the algebra $\approxexists{\frk{A}}(B') = (B', \lozenge, \exists')$.
\end{definition}


\begin{lemma}
  \label{lem:approximate-e-is-ms4-algebra}
  With $\frk{A}, B'$ as above, $\approxexists{\frk{A}}(B')$ is an $\msf{MS4}$-algebra.
\end{lemma}

\begin{proof}
  (1) It is well-known that $\exists'$ is an $\msf{S4}$-operator on $B'$ (\cite[Lem.~2.3]{MT44}).
  To show it is an $\msf{S5}$-operator, simply note that if $a \in B'$, $\forall' a \in B^0 \cap B'$ by definition, so $\exists' \forall' a = \exists \forall' a = \forall' a$.

  We now show $\exists' \lozenge \exists' a = \lozenge \exists' a$ which is equivalent to the $\msf{MS4}$-axiom (see, e.g., \cite[Lem.~2.4]{BM23}).
  $\exists' a$ is a finite meet of $\exists$-fixpoints and hence an $\exists$-fixpoint.
  Then since $\exists \lozenge \exists = \lozenge \exists$ holds in $\frk{A}$, $\lozenge \exists' a$ is also an $\exists$-fixpoint, and belongs to $B'$ since $B'$ is a $\lozenge$-subalgebra.
  Then $\lozenge \exists' a \in B^0 \cap B'$, which means $\exists' \lozenge \exists' a = \exists \lozenge \exists' a = \lozenge \exists' a$.
\end{proof}

\begin{theorem}
  $\bb{MS4}[n]$ has the fmp.
\end{theorem}

\begin{proof}
  If $\bb{MS4}[n] \not\models \varphi$, then there is an $\msf{MS4}[n]$-algebra $\frk{A}$ and an $n$-tuple $\bar{a}$ so that $\varphi(\bar{a}) \neq 1$ in $\frk{A}$.
  Let $S$ be the set of subterms of $\varphi$.
  Let $B' = \la S \ra_{(B, \lozenge)}$ be the subalgebra of the $\lozenge$-reduct of $\frk{A}$ generated by $S$.
  $B'$ is finite since the variety $\bb{S4}[n]$ is locally finite.
  By \cref{lem:approximate-e-is-ms4-algebra} $\approxexists{\frk{A}}(B')$ is a finite $\msf{MS4}$-algebra and, since the identities $P_n$ involve only $\lozenge$, $\approxexists{\frk{A}}(B') \models P_n$ since its $\msf{S4}$-reduct is a subalgebra of the $\msf{S4}$-reduct of $\frk{A}$.
  Clearly for any $a$ with $\exists a \in B'$, $\exists' a = \exists a$; hence the computation of $\varphi(\bar{a})$ in $\approxexists{\frk{A}}(B')$ is identical to that in $\frk{A}$ and in particular $\varphi(\bar{a}) \neq 1$ in $\approxexists{\frk{A}}(B')$.
  Then we have identified a finite $\msf{MS4}[n]$-algebra in which $\varphi(\bar{a}) \neq 1$.
\end{proof}

\begin{theorem}
  $\bb{MS4_S}[n]$ has the fmp.
\end{theorem}

\begin{proof}
  Recall that $\bb{MS4_S}[n] = \bb{MS4} + P^0_1 + P_n$.
  If $\bb{MS4_S}[n] \not\models \varphi$, then there is an $\msf{MS4_S}[n]$-algebra $\frk{A}$ and an $n$-tuple $\bar{a}$ so that $\varphi(\bar{a}) \neq 1$ in $\frk{A}$.
  Let $S$ be the set of subterms of $\varphi$.
  Let $B' = \la S \ra_{(B, \lozenge, \blacklozenge)}$ be the subalgebra of $(B, \lozenge, \blacklozenge)$ generated by $S$.
  $\blacklozenge$ is an $\msf{S5}$-operator (since $P^0_1$ holds in $\frk{A}$).
  Moreover, the identity $\lozenge a \leq \blacklozenge a$ (i.e., $\lozenge a \leq \lozenge \exists a$) holds in $\frk{A}$ (since $\exists$ and $\lozenge$ are increasing).
  Then the algebra $(B, \lozenge, \blacklozenge)$ is an $\msf{MS4}$-algebra belonging to the variety $\bb{S4}_u$ ($\bb{S4}$ with a universal modality, see \cite{GVP92,Ben96} and \cite[Def.~3.10]{BM23}).
  This variety is locally finite by \cite[Cor.~4.5]{BM23}, and hence $B'$ is finite.
  Since the identities $P_n$ and $P^0_1$ involve only $\lozenge, \blacklozenge$, we have $\approxexists{\frk{A}}(B') \models P_n, P^0_1$, since $(B', \lozenge, \blacklozenge)$ is an $\msf{MS4}$-subalgebra of $(B, \lozenge, \blacklozenge)$.
  By the same reasoning as the previous theorem, the computation of $\varphi(\bar{a})$ in $\approxexists{\frk{A}}(B')$ is identical to that in $\frk{A}$, thus we have a finite $\msf{MS4_S}[n]$-algebra in which $\varphi(\bar{a}) \neq 1$.
\end{proof}

\section{The Casari axiom and \texorpdfstring{$\mpsfl$}{M+S4}}
\label{sec:casari}

\subsection{The dual spaces of \texorpdfstring{$\mpsfl$}{M+S4}-algebras}
\label{sec:casari-dual-spaces}

Casari's formula is the following first-order predicate formula:
\[
  \msf{Cas} = \forall x \, ((P(x) \to \forall y \, P(y)) \to \forall y \, P(y)) \to \forall x \, P(x)
\]
In \cite{BBI23} it was shown that the monadic version of Casari's formula
\[
  \msf{MCas} = \forall ((p \to \forall p) \to \forall p) \to \forall p
\]
is necessary to obtain a faithful interpretation of monadic intuitionistic predicate logic ($\msf{MIPC}$) into monadic G\"odel--L\"ob logic ($\msf{MGL}$) and hence a faithful provability interpretation of $\msf{MIPC}$.
In light of this, it is natural to study the logic obtained by asserting the G\"odel translation of this formula (see, e.g., \cite[Def.~4.5,Rem.~4.6]{BBI23}).

\begin{definition}
  Let
  \[
    \msf{M^{+}Cas} = \blacksquare(\square(\square p \to \blacksquare p) \to \blacksquare p) \to \blacksquare p
  \]
  and define the subvariety $\mpsfv = \bb{MS4} + \msf{M^{+}Cas}$.
\end{definition}

We now move to establish some fundamental properties of descriptive $\mpsfl$-frames.
\cref{thm:cas-frame-condition} generalizes the analogous result obtained in \cite[Lem.~4.8]{BBI23} for $\mpgrzl$-frames.

\begin{definition}
  Let $\frk{F} = (X, R, E)$ be a descriptive $\msf{MS4}$-frame and $A \subseteq X$.
  We call $A$ \textit{flat} if, for all $x, y \in A$, $x R y$ implies $y R x$ (i.e., $R \vert_A$ is an equivalence).
\end{definition}

We will refer to an equivalence class $E(x)$ as an \textit{$E$-cluster}.
We call a subset $A \subseteq X$ \textit{$E$-saturated} if it is an $E$-upset ($x \in A \tand x E y \to y \in A$) or, equivalently, a union of $E$-clusters.

\begin{lemma}
  \label{lem:flat-cluster-lemmas}
  Let $\frk{F} = (X, R, E)$ be an $\msf{MS4}$-frame, and suppose $A \subseteq X$ is an $E$-saturated set.
  \begin{enumerate}
    [label=\normalfont(\arabic*), ref = \thelemma(\arabic*)]
    \item \label[lemma]{lem:flat-cluster-lemmas:1} If $m \in \qmax A$ and $E(m)$ is flat, then $x \in \qmax A$ for all $x \in E(m)$.
    \item \label[lemma]{lem:flat-cluster-lemmas:2} If $m \in \qmax A$, then for all $x \in A$, $mQx$ implies $xQm$ (i.e., $m$ is also $Q$-quasimaximal in $A$).
  \end{enumerate}
\end{lemma}

\begin{proof}
  (1)
  Let $x \in E(m)$ and $x R y$ for $y \in A$.
  Then $m E x R y$, so $m R y' E y$ for some $y'$ and since $A$ is $E$-saturated, $y' \in A$. 
  Then since $m \in \qmax A$, $y' R m$.
  So we have $y E y' R m$ and hence $y R x' E m$ for some $x'$, and $x' \in E(m)$.
  Then we have $x R y R x'$ and, since $E(m)$ is flat, $x' R x$.
  We conclude that $y R x' R x$ and $x \in \qmax A$.

  (2)
  Let $x \in A$ and $m Q x$.
  Then $m R x' E x$ and since $A$ is $E$-saturated, $x' \in A$.
  Since $m \in \qmax A$, $x' R m$.
  Then $x E x' R m$ and so $x R m' E m$ for some $m'$, and $x Q m$.
\end{proof}

The following characterizes the (descriptive) frame condition for the formula $\msf{M^{+}Cas}$, which generalizes the one obtained for $\msf{M^{+}Grz}$-frames in \cite[Lem.~4.8]{BBI23}.

\begin{theorem}
  \label{thm:cas-frame-condition}
  Let $\frk{F}$ be an $\msf{MS4}$-frame.
  $\frk{F} \models \msf{M^{+}Cas}$ iff, for any $E$-saturated clopen $U$, $x \in \qmax U$ implies $E(x)$ is flat.
\end{theorem}

\begin{proof}
  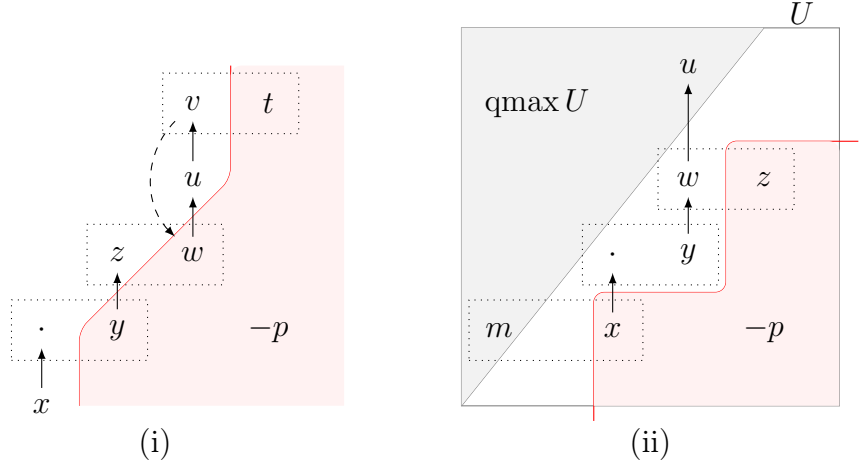
\begin{figure}[h]
    \centering
    \begin{tikzpicture}[
      scale=1,
      dot/.style={circle,fill=black,minimum size=6pt,inner sep=0}
    ]
      \draw[thin,red,rounded corners] (0.5, 0) -- (0.5, 1) -- (2.5, 3) -- (2.5, 4.5);
      \fill[red!5] [rounded corners] (0.5, 0) -- (0.5, 1) -- (2.5, 3) -- (2.5, 4.5) [sharp corners] -- (4, 4.5) -- (4, 0) -- cycle;

      \node (x) at (0,0) {$x$};
      \node (xy) at (0,1) {$\cdot$};
      \node (y) at (1,1) {$y$};
      \node (z) at (1,2) {$z$};
      \node (w) at (2,2) {$w$};
      \node (u) at (2,3) {$u$};
      \node (v) at (2,4) {$v$};
      \node (t) at (3,4) {$t$};

      \node at (3, 1) {$-p$};

      \draw[-latex] (x) -- (xy);
      \draw[-latex] (y) -- (z);
      \draw[-latex] (w) -- (u);
      \draw[-latex] (u) -- (v);
      \draw[-latex,dashed] (v) to[out=225,in=135] (w);

      \draw[dotted] (-0.4, 0.6) rectangle (1.4,1.4);
      \draw[dotted] (0.6, 1.6) rectangle (2.4,2.4);
      \draw[dotted] (1.6, 3.6) rectangle (3.4,4.4);  
      
      \node at (1.5, -0.5) {(i)};
    \end{tikzpicture}
    \hspace{0.5in}
    \begin{tikzpicture}[
      scale=1,
      dot/.style={circle,fill=black,minimum size=6pt,inner sep=0}
    ]
      \draw[thin,gray,rounded corners] (-1,-1) -- (3,4);
      \draw[thin,gray] (-1,-1) rectangle (4, 4);
      \fill[black!5] (-1, -1) -- (-1, 4) -- (3, 4) -- cycle;
      \draw[thin,red,rounded corners] (0.75,-1.2) -- (0.75, 0.5) -- (2.5, 0.5) -- (2.5, 2.5) -- (4.25, 2.5);
      \fill[red!5] [rounded corners] (0.75, -1) -- (0.75, 0.5) -- (2.5, 0.5) -- (2.5, 2.5) -- (4, 2.5) [sharp corners] -- (4, -1) -- cycle;

      \node (m) at (-0.5,0) {$m$};
      \node (x) at (1,0) {$x$};
      \node (xy) at (1,1) {$\cdot$};
      \node (y) at (2,1) {$y$};
      \node (w) at (2,2) {$w$};
      \node (z) at (3,2) {$z$};
      \node (u) at (2,3.5) {$u$};

      \node at (3, 0) {$-p$};
      \node at (0, 3) {$\qmax U$};
      \node at (3.5, 4.2) {$U$};

      \draw[-latex] (x) -- (xy);
      \draw[-latex] (y) -- (w);
      \draw[-latex] (w) -- (u);

      \draw[dotted] (-0.9, -0.4) rectangle (1.4,0.4);
      \draw[dotted] (0.6, 0.6) rectangle (2.4,1.4);
      \draw[dotted] (1.6, 1.6) rectangle (3.4,2.4);

      \node at (1.5, -1.5) {(ii)};
    \end{tikzpicture}
    \caption{
      The arguments for both direction of \cref{thm:cas-frame-condition}.
    }
    \label{fig:cas-frame-condition}
  \end{figure}

  Suppose $\frk{F} \not\models \msf{M^{+}Cas}$ (see \cref{fig:cas-frame-condition} (i)).
  Then take a clopen valuation $\nu$ and $x \in X$ such that $x \not\models \msf{M^{+}Cas}$.
  That is, $x \not\models \blacksquare p$ but $x \models \blacksquare (\square (\square p \to \blacksquare p) \to \blacksquare p)$.
  Then we may find $x Q y$ with $y \not\models p$ and $y R z$ with $z \models \qmax E(\bbrackets{-p})$ (by \cref{thm:sees-maximal:1}).
  We will show that $E(z)$ is non-flat, falsifying the frame condition.
  Suppose toward contradiction that $E(z)$ is flat.
  Since $z \in E(\bbrackets{-p})$, we can find $z E w$ for $w \models -p$.
  Since $x Q w$ and $w \not\models \blacksquare p$, we must have $w \not\models \square(\square p \to \blacksquare p)$.
  Then we may find $w R u$ for $u \models \square p$ but $u \not\models \blacksquare p$; then we may find $u R v E t$ with $t \not\models p$; note that $v \in E(\bbrackets{-p})$.
  Since $z \models \qmax E(\bbrackets{-p})$ and $E(z)$ is flat, by \cref{lem:flat-cluster-lemmas:1} we have $w \models \qmax E(\bbrackets{-p})$.
  Then since $w R v$, we have $v R w$.
  But then $u R w$, which is a contradiction since $u \models \square p$, but $w \not\models p$.

  Suppose $m \in \qmax U$ for an $E$-saturated clopen $U$, and $E(m)$ is non-flat (see \cref{fig:cas-frame-condition} (ii)).
  Then we may find $x \in E(m)$ (and $x \in U$) with $x R x'$ and $\neg(x' R x)$; in particular, $x \not\in \qmax U$.
  Since $\qmax U$ is closed, we may separate $x$ from $\qmax U$ with a clopen set $V$.
  Define a valuation by $\nu(p) = -(U \cap V)$; note that $x \not\models p$ and $\qmax U \subseteq \bbrackets{p}$.
  We now claim that $x \not\models \msf{M^{+}Cas}$ (and hence $\frk{F} \not\models \msf{M^{+}Cas}$).
  We have that $x \not\models \blacksquare p$, so we are required to show that $x \models \blacksquare(\square(\square p \to \blacksquare p) \to \blacksquare p)$.
  That is, $x Q y$ implies $y \models \square(\square p \to \blacksquare p) \to \blacksquare p$.
  So suppose $x Q y$ and $y \not\models \blacksquare p$.
  Then we may find $y R w E z$ with $z \not\models p$.
  Then we have $z \in U$ (since $\bbrackets{-p} = U \cap V$) and $w \in U$ (since $U$ is $E$-saturated).
  Then we may find $w R u$ for $u \in \qmax U$ again by \cref{thm:sees-maximal:1}; note $u \models p$.
  We claim that $u \models \square p$ but $u \not\models \blacksquare p$ (which implies, since $y R u$, that $y \not\models \square(\square p \to \blacksquare p)$).
  Indeed, to show $u \models \square p$, suppose $u R u'$ with $u' \not\models p$. Then $u' \in U$ and, since $u \in \qmax U$, $u' R u$ and $u' \in \qmax U$, but this is a contradiction since $\qmax U$ and $\bbrackets{-p}$ are disjoint.
  On the other hand, we have $m \in \qmax U, u \in U$, and $m Q u$.
  By \cref{lem:flat-cluster-lemmas:2}, $u Q m$ and hence $u Q x$; since $x \not\models p$, $u \not\models \blacksquare p$, completing the claim.
\end{proof}

We will now establish several properties concerning the layers of $\mpsfl$-frames (introduced in \cref{sec:prelims-depth}).
For an arbitrary $\msf{MS4}$-frame $\frk{F}$, \cref{thm:sees-maximal:2} says that the first layer $D_1(\frk{F})$ is non-empty and closed; if $\frk{F}$ is infinite depth, it may very well be the case that some or all of the successive layers are empty.
However, we will have reason to speak of the finite layers of $\frk{F}$ that \textit{are} non-empty; for this purpose, we let $\ell(\frk{F})$ denote the greatest integer $n \geq 1$ such that $D_n(\frk{F}) \neq \varnothing$ or, if there is no such $n$, $\ell(\frk{F}) = \omega$.
Certainly $\ell(\frk{F}) \leq d(\frk{F})$ and, if $\frk{F}$ is of finite depth, $\ell(\frk{F}) = d(\frk{F})$; however, we may have $d(\frk{F}) = \omega$ while $\ell(\frk{F}) < \omega$, for example in the $\msf{S4}$-frame $(\omega + 1, \leq)$.

We first note that successive unions of non-empty layers of an $\msf{MS4}$-frame are closed; this was established for the dual spaces of monadic Heyting algebras in \cite[Lemma 7]{Bez00}; here we adapt the result to $\msf{S4}$-frames and hence $\msf{MS4}$-frames:

\begin{lemma}
  \label{lem:layers-closed}
  Let $\frk{F} = (X, R, E)$ be a descriptive $\msf{S4}$-frame.
  Then the sets $\bigcup_{i = 0}^n D_i(\frk{F})$ are closed for any $n < \ell(\frk{F})$.
\end{lemma}

\begin{proof}
  Clearly $D_0(\frk{F}) = \varnothing$ is closed and the case $n=1$ is \cref{thm:sees-maximal:2}.
  Suppose $M = \bigcup_{i=1}^{n-1} D_i(\frk{F})$ for $n \geq 1$ is closed.
  To show that $N = \bigcup_{i=1}^n D_i(\frk{F})$ is closed, take a point $x \in -N$.
  By definition, there must be some $y \in D_n(\frk{F})$ such that $x \, R^+ \, y$ (since $D_n(\frk{F})$ is non-empty).
  Take $y \in U$ for a clopen neighborhood $U$ disjoint from the closed set $M$.
  Since $\neg(y \, R \, x)$, we may find a clopen $R$-upset $V$ with $y \in V$ and $x \not\in V$ (\cite[Thm.~3.1.2(IV)]{Esa19}).
  Then $W = R^{-1}(U \cap V) - V$ is a clopen neighborhood of $x$; we claim it is disjoint from $N$.
  By construction, $x$ belongs to this set; suppose $z \in N \cap W$.
  Then $z \, R \, w$ for $w \in U \cap V$, and hence $d(z) \geq d(w)$.
  But $d(z) \leq n$, while $d(w) \geq n$ (since $U \cap V$ is disjoint from $M$), so the only possibility is that $z, w \in D_n(\frk{F})$, and $w \, R \, z$. 
  But $w \in V$ and $V$ was an $R$-upset, meaning that $z \in V$, contradictory to $z \in W$.
  Then $W \cap N = \varnothing$, and $N$ is closed.
\end{proof}

The layers of $\mpsfl$-frames are vastly better behaved than those of arbitrary $\msf{MS4}$-frames, as the following theorems establish.

\begin{lemma}
  \label{lem:layers-saturated}
  Suppose $\frk{F} = (X, R, E)$ be an $\mpsfl$-frame.
  Then $D_n(\frk{F})$ is $E$-saturated, for each $n < \ell(\frk{F})$.
\end{lemma}

\begin{proof}
  By induction on $n$; when $n=0$, we have $D_0(\frk{F}) = \varnothing$ which is $E$-saturated.
  Suppose that each $D_i(\frk{F})$ for $i < n$ is $E$-saturated and take $x \in D_n(\frk{F})$.
  By the previous lemma, $S = \bigcup_{i=0}^{n-1} D_i(\frk{F})$ is a closed set, and it is $E$-saturated by assumption (as a union of $E$-saturated sets).
  Then $-S$ is an open $E$-saturated set.
  Hence, there is a clopen $U \subseteq -S$ with $x \in U$ and furthermore $E(U)$ is an $E$-saturated clopen contained in $-S$.
  Next, observe that $x \in \qmax E(U)$, since $x R^+ y$ implies $y \in S$.
  By \cref{thm:cas-frame-condition}, $E(x)$ is flat.
  Finally, to show $D_n(\frk{F})$ is $E$-saturated, suppose $x E y$; we conclude by \cref{lem:flat-cluster-lemmas:1} that $y \in \qmax E(U)$, and hence $y \in D_n(\frk{F})$.
\end{proof}

\begin{lemma}
  \label{lem:fin-depth-cluster-flat}
  Let $\frk{F} = (X, R, E)$ be an $\mpsfl$-frame.
  If $x$ is a point of finite depth, then $E(x)$ is flat.
\end{lemma}

\begin{proof}
  We have $x \in D_n(\frk{F})$ for some $n$ which is $E$-saturated.
  So for all $y \in E(x)$, $y \in D_n(\frk{F})$ and hence $x R y$ implies $y R x$.
\end{proof}

\begin{definition}
  Suppose $\frk{F} = (X, R, E)$ is a descriptive $\msf{MS4}$-frame and $U \subseteq X$ a clopen set.
  Define
  \begin{itemize}
    \item $R \vert_U = R \cap U^2$
    \item $E \vert_U = E \cap U^2$
    \item $\frk{F} \vert_U = (U, R\vert_U, E \vert_U)$
  \end{itemize}
  We will be particularly interested in the case when $U = D_n(\frk{F})$ (and $D_n(\frk{F})$ is clopen) for $n \leq \ell(\frk{F})$; in this case, we write $\frk{F} \vert_n$ to mean $\frk{F} \vert_{D_n(\frk{F})}$.
  We also apply the same notation to algebras; if $\frk{A}$ is an $\msf{MS4}$-algebra, $\frk{A} \vert_U$ denotes $(\frk{F} \vert_U)^*$.
\end{definition}

\begin{remark}
  \label{rem:relativization-clopen}
  When $U$ is clopen, the algebra $\frk{A} \vert_U$ is the \textit{relativization} to the element $U$ -- that is, the algebra
  \[
    ([0, U], \cap, \cup, (-)', \varnothing, U, \lozenge', \exists')
  \]
  where the relative complementation and modal operators are given by
  \[
    a' = U - a \quad \lozenge' a = \lozenge a \wedge U \quad \exists' a = \exists a \wedge U
  \]
  and in this case, $\frk{F}_U$ is the dual space of $\frk{A}_U$.
\end{remark}

Recall that $\bb{S5}^2 = \bb{MS4}[1]$ is the variety of Boolean algebras with two commuting $\msf{S5}$-operators (or, dually, descriptive frames with two commuting equivalence relations).
The algebras $\frk{A}_n$ for a general $\msf{MS4}$-frame will in general only belong to $\bb{S5}_2$, the variety of Boolean algebras with possibly non-commuting $\msf{S5}$-operators (or, dually, descriptive frames with two possibly non-commuting equivalence relations) (see \cref{sec:translation} or \cite[Sec.~6]{BM23}).
For $\mpsfl$-frames however, the situation is much better:

\begin{lemma}
  \label{lem:layers-commute}
  Let $\frk{F} = (X, R, E)$ be an $\mpsfl$-frame, and $n \leq d(\frk{F})$ such that $D_n(\frk{F})$ is clopen.
  Then $\frk{F} \vert_n$ is a (descriptive) $\msf{S5}^2$-frame.
\end{lemma}

\begin{proof}
  Since $D_n(\frk{F})$ is clopen subset of a Stone space, it is also a Stone space in the subspace topology, and the restrictions $R_n$ and $E_n$ of the relations are continuous on $D_n(\frk{F})$ (as they are defined by the intersection with a clopen set).
  Suppose $x \, E \vert_n \, y \, R \vert_n \, z$, so $x,y,z \in D_n(\frk{F})$.
  Then $x R y' E z$ for some $y' \in X$, but since $D_n(\frk{F})$ is $E$-saturated, we have $y' \in D_n(\frk{F})$.
  Hence $x \, R \vert_n \, y' \, E \vert_n \, z$.
  Then $R \vert_n \circ E \vert_n \subseteq E \vert_n \circ R \vert_n$ and since both are equivalence relations, they commute.
\end{proof}

The following lemma expresses, in several ways, the fact that the $R$-depth and $Q$-depth of $\mpsfl$-frames (and varieties thereof) coincide (see \cref{sec:prelims-depth} for the concepts and associated notation used in the following results).

\begin{lemma}
  \label{lem:depths-coincide-points}
  Let $\frk{F} = (X, R, E)$ be an $\mpsfl$-frame, and $x \in \frk{F}$ a point of finite depth.
  Then 
  \begin{enumerate}
    [label=\normalfont(\arabic*), ref = \thelemma(\arabic*)]
    \item \label[lemma]{lem:depths-coincide-points:1} $x R^+ y$ implies $E(x) \, \overline{R}^+ \, E(y)$.
    \item \label[lemma]{lem:depths-coincide-points:2} $d^0(x) = d(x)$
  \end{enumerate}
\end{lemma}

\begin{proof}
  (1) If $x R^+ y$, we have $E(x) \, \overline{R} \, E(y)$ but we cannot have $E(y) \, \overline{R} \, E(x)$; if so, then $y R x'$ for some $x' \in E(x)$, and hence $x R^+ x'$, but by \cref{lem:fin-depth-cluster-flat}, $E(x)$ is flat.
  
  (2) It is always the case that $d^0(x) \leq d(x)$, since any choice of representatives of an $\overline{R}^+$-chain in $\frk{F}^0$ form an $R^+$-chain in $\frk{F}$ (since $RE \subseteq ER$).
  Conversely, by (1), any $R^+$-chain emanating from $x$ yields an $\overline{R}^+$-chain in $\frk{F}^0$, so $d(x) \leq d^0(x)$.
\end{proof}

This result is enough to imply that if $d(\frk{F}) < \omega$, then $d^0(\frk{F}) = d(\frk{F})$.
In fact, $d(\frk{F})$ and $d^0(\frk{F})$ coincide for arbitrary $\mpsfl$-frames, but the case $d(\frk{F}) = \omega$ requires more care:

\begin{lemma}
  \label{lem:depths-coincide-frames}
  Let $\frk{F} = (X, R, E)$ be an $\mpsfl$-frame.
  Then $d^0(\frk{F}) = d(\frk{F})$.
\end{lemma}

\begin{proof}
  If $d(\frk{F})$ is finite, then $d^0(x) = d(x)$ for all $x$ in $\frk{F}$ by \cref{lem:depths-coincide-points}, so $d^0(\frk{F}) = d(\frk{F})$.
  
  Now suppose $d(\frk{F}) = \omega$, and that there exist points in $\frk{F}$ of arbitrarily large finite depth.
  Then by \cref{lem:depths-coincide-points}, there exist points of arbitrarily large $Q$-depth, and $d^0(\frk{F}) = \omega$ as well.
  
  Otherwise, suppose $d(\frk{F}) = \omega$ but $\frk{F}$ does not contain points of arbitrarily high finite depth; suppose the depth of any (finite-depth) point in $\frk{F}$ is at most $m$.
  Then there must exist points of infinite depth in $\frk{F}$.
  Suppose toward contradiction that $d^0(\frk{F}) < \omega$.
  Let $n$ be the smallest integer such that there exists a point $x$ with $d(x) = \omega$ and $d^0(x) = n$ (that is, any point of strictly smaller $Q$-depth than $x$ has finite $R$-depth of at most $m$).
  We will establish two claims:
  \begin{enumerate}
    \item There exists $y \in E(x)$ with $x \, R^+ \, y$ (i.e., $E(x)$ is not flat)
    \item There exists $z \in E(x)$ such that $z \in \qmax E(U)$ for a clopen set $U$.
  \end{enumerate}
  Together, this contradicts the frame condition for $\msf{MSCas}$ in \cref{thm:cas-frame-condition}, since $z$ is quasimaximal in an $E$-saturated clopen, but $E(z) = E(x)$ is not flat; hence this case cannot occur.

  (1) It must be the case that $x \, R^+ \, x'$ for some $x'$ with $d^0(x') = d^0(x)$ -- otherwise, $x \, R^+ \, x'$ implies $d^0(x') < d^0(x)$ and hence, by the assumption of minimality, that $d(x') \leq m$; but then $d(x) \leq m+1$, contradictory to the choice of $x$.
  Then, since $x \, R^+ \, x'$, $x \, Q \, x'$, and we must have $x' \, Q \, x$, i.e. $x' \, R \, y$ for $y \in E(x)$. 
  Evidently, $x \, R^+ \, y$.

  (2) Since $E(x)$ is closed, we may find $x \, R \, z$ for $z \in \qmax E(x)$, by \cref{thm:sees-maximal:1}.
  Clearly $d^0(z) = d^0(x) = n$.
  As $\frk{F}^0$ is an $\msf{S4}$-frame of finite depth, the set $M^0 = \bigcup_{i=1}^{n-1} D_i(\frk{F}^0)$ is closed in $X/E$, by \cref{lem:layers-closed}.
  Letting $\pi : X \to X/E$ be the quotient map, $M = \pi^{-1}(M^0)$ is closed in $X$ (simply by continuity), and consists precisely of the points of $Q$-depth strictly less than $n$.
  Then we may take $y \in U$ for a clopen set $U$ disjoint from $M$.
  Since $M$ is $E$-saturated (being a $Q$-upset), $E(U)$ is an $E$-saturated clopen set disjoint from $M$.
  Evidently $z \in \qmax_Q E(U)$ (since $E(U)$ is disjoint from $M$, it cannot contain points of lower $Q$-depth).
  We claim it actually belongs to $\qmax E(U)$.
  Suppose $z \, R \, w$ for $w \in E(U)$.
  Then $z \, Q \, w$, and hence $w \, Q \, z$; i.e. $w \, R \, z'$ for $z' \in E(z)$.
  But $E(z) = E(x)$, and $z \in \qmax E(x)$, so $z' \, R \, z$.
  It follows that $w \, R \, z$, and $z \in \qmax E(U)$.
\end{proof}

As an immediate consequence, we obtain:

\begin{theorem}
  \label{thm:depths-coincide-varieties}
  Let $\bb{V} \subseteq \mpsfv$.
  Then $d(\bb{V}) \leq n$ iff $d^0(\bb{V}) \leq n$; equivalently, $\bb{V} \models P_n$ iff $\bb{V} \models P^0_n$.
\end{theorem}

The following lemma states that in an s.i. $\mpsfl$-frame $\frk{F}$ of finite depth, the set of $Q$-roots of $\frk{F}$ coincides with the bottom layer of $\frk{F}$, which will be clopen.

\begin{lemma}
  \label{lem:si-bottom-clopen}
  Let $\frk{F}$ be a s.i.~$\mpsfl$-frame with $d(\frk{F}) = h$.
  Then $T(\frk{F}) = D_h(\frk{F})$ and $D_h(\frk{F})$ is clopen.
\end{lemma}

\begin{proof}
  $\frk{F}$ is strongly $Q$-rooted, so there is an open set of $Q$-roots and in particular there exists a $Q$-root $r$.
  Since $\bigcup_{i=1}^{h-1} D_i(\frk{F})$ is a closed $Q$-upset (since it is $E$-saturated, by \cref{lem:layers-saturated}), we must have $r \in D_h(\frk{F})$ (otherwise, $Q(r) \subset X$).
  Now suppose $x \in D_h(\frk{F})$.
  We have $rQx$ and hence $r R z E x$ for some $z \in X$, so $z \in D_h(\frk{F})$ since $D_h(\frk{F})$ is $E$-saturated.
  Then $r \, R \vert_h \, z \, E \vert_h \, x$ and hence $r \, E \vert_h \, z' \, R \vert_h \, x$ since $R \vert_h$ and $E \vert_h$ commute.
  This means that $x Q r$ and $x$ is also a $Q$-root of $\frk{F}$.
  
  This shows that the set of $Q$-roots of $\frk{F}$ is precisely $D_h(\frk{F})$.
  We then note that an open set of roots is actually clopen, by the following standard argument: Since $T(\frk{F})$ is open and nonempty, there is a clopen $U \subseteq T(\frk{F})$, and $T(\frk{F}) = Q^{-1}(U)$ which is clopen since $Q$ is a continuous quasiorder.
\end{proof}

\subsection{A semantic characterization of local finiteness}

In light of \cref{lem:si-bottom-clopen}, we make the following definition:

\begin{definition}
  Let $\frk{A}$ be a s.i.~$\mpsfl$-algebra of finite depth.
  \begin{enumerate}
    \item Let $T(\frk{A})$ denote the (clopen) element of $\frk{A}$ corresponding to $T(\frk{A}_*)$.
    \item Let $\frk{A} \vert_T$ (resp. $\frk{F} \vert_T$ for frames) denote $\frk{A} \vert_{d(\frk{A})}$ (resp. $\frk{F} \vert_{d(\frk{F})}$).
  \end{enumerate}
\end{definition}

\begin{definition}
  Let $\bb{V} \subseteq \mpsfv$ be a variety of $\mpsfl$-algebras of finite depth.
  Define
  \[
    \bb{V} \vert_T = \bb{Var}(\setbuilder{\frk{A} \vert_T}{\frk{A} \in \bb{V}_{\text{SI}}}) \subseteq \bb{S5}^2
  \]
\end{definition}

We think of $\bb{V}_T$ as the ``variety of bottom layers of s.i. $\bb{V}$-algebras'' (note that $\bb{V}_T \subseteq \bb{S5}^2$ in light of \cref{lem:layers-commute}).
This construction allows us to give a first characterization of local finiteness in $\mpsfv$.

\begin{theorem}
  \label{thm:mpsfl-semantic-criterion}
  $\bb{V} \subseteq \mpsfv$ is locally finite iff
  \begin{enumerate}
    \item $\bb{V}$ is finite depth
    \item $\bb{V} \vert_T$ is locally finite.
  \end{enumerate}
\end{theorem}

\begin{proof}
  Suppose $\bb{V}$ is locally finite, so that we have $f : \omega \to \omega$ a bound on the size of an $m$-generated algebra from $\bb{V}$.
  Since the $\bb{S4}$-reduct of $\bb{V}$ is a uniformly locally finite class (\cite[Def.~2.1(c)]{Bez01}), it generates a locally finite variety (by \cite[Thm.~3.7(2)]{Bez01}) and $\bb{V}$ is finite depth, so (1) holds.
  To establish (2), we show that the given generating class of $\bb{V} \vert_T$ (which is well-defined, as $\bb{V}$ is finite depth) is uniformly locally finite.
  Indeed, let $\frk{A} \in \bb{V}_\text{SI}$ with dual frame $\frk{F}$ and $\frk{B} = \la g_1, \dots, g_m \ra_{\frk{A} \vert_T}$ be an $m$-generated subalgebra of $\frk{A} \vert_T$. 
  The relative operations of $\frk{A} \vert_T$ are definable in $\frk{A}$ using the element $T(\frk{A})$ as explained in \cref{rem:relativization-clopen}.
  Hence the subalgebra $\la g_1, \dots, g_m, T(\frk{A}) \ra_\frk{A}$ of $\frk{A}$ contains, in particular, every element of $\frk{B}$ and has size bounded by $f(m + 1)$.
  Then by \cite[Thm.~3.7(2)]{Bez01} $\bb{V}_T$ is locally finite.

  Suppose (1)-(2) hold.
  Let $f : \omega \to \omega$ be a uniform bound on the size of an $n$-generated algebra from $\bb{V} \vert_T$, and $d(\bb{V}) = h$.
  We will show that $\bb{V} + P^0_n$ is locally finite for each $1 \leq n \leq h$, and then note that $\bb{V} + P^0_h = \bb{V}$ (since $d^0(\bb{V}) \leq d(\bb{V}) = h$).

  In the base case, we have $\bb{V} + P^0_1 \subseteq \bb{S5}^2$ (since by \cref{lem:depths-coincide-frames}, $\bb{V} \models P_1$). 
  Moreover, $\bb{V} + P^0_1 \subseteq \bb{V} \vert_T$ since every s.i.~algebra $\frk{A}$ from $\bb{V} + P^0_1$ has $\frk{A} \vert_T = \frk{A}$ and hence $\frk{A}$ is among the generators of $\bb{V} \vert_T$.
  Then by the assumption (2), $\bb{V} + P^0_1$ is locally finite.

  Suppose $\bb{V} + P^0_{n-1}$ is locally finite so that we have $g : \omega \to \omega$ a bound on the size of an $m$-generated algebra from $\bb{V} + P^0_{n-1}$.
  Let $\frk{A}$ be an s.i.~algebra from $\bb{V} + P^0_n$ generated by $G = \set{g_1, \dots, g_m}$, and $\frk{F}$ its dual space.
  Set $T = T(\frk{A})$.
  We have $\frk{A} \cong \frk{A} \vert_T \times \frk{A}'$ as a Boolean algebra, where $\frk{A}'$ is the Boolean relativization to the element $-T$ (identifying $\frk{A}$ with the algebra of clopen subsets of $\frk{F}$, $\frk{A} \vert_T$ consists of clopen subsets of $T$ and $\frk{A}'$ of clopen subsets of $-T$).
  Indeed $-T$ is the largest proper $Q$-upset of $\frk{F}$ (since $T$ is the set of $Q$-roots of $\frk{F}$) and hence $\frk{A}'$ is a homomorphic image of $\frk{A}$ from $\bb{V} + P^0_{n-1}$ also generated by $G$ (since closed $Q$-upsets of $\frk{F}$ are in correspondence with modal filters and congruences of $\frk{A}$; see, e.g., \cite[Thm.~3.4]{BM23}).
  By the induction hypothesis, $\abs{\frk{A}'} \leq g(m)$.
  We claim that $\frk{A} \vert_T$ is generated by the set $G \vert_T \cup G'$, where $G \vert_T = \set{g_1 \cap T, \dots, g_m \cap T}$ and $G' = \setbuilder{R^{-1}(U) \cap T}{U \in \frk{A}'}$ (here and in the rest of the proof, we freely identify the algebra $\frk{A} \vert_T$ with clopen subsets of $T$, and $\frk{A}'$ with clopen subsets of $-T$).

  Indeed every element of $\frk{A} \vert_T$ is generated in $\frk{A}$ by $G$, so it suffices to show that for every $A \in \frk{A}$, $A \cap T \in \la G \vert_T \cup G' \ra_{\frk{A} \vert_T}$ or, for every $\msf{MS4}$-polynomial $p(x_1, \dots, x_m)$, we have $p^\frk{A}(g_1, \dots, g_m) \cap T \in \la G \vert_T \cup G' \ra_{\frk{A} \vert_T}$.
  We do this by induction on $p$.

  Suppose $p = x_i$ (or $p = 0$ or $p = 1$).
  Then $p^\frk{A}(\bar{g}) \cap T = g_i \cap T \in G \vert_T$.

  Suppose $p = q \wedge r$ where $q$ and $r$ satisfy the induction hypothesis.
  Then,
  \[
    p^\frk{A}(\bar{g}) \cap T = (q^\frk{A}(\bar{g}) \cap r^\frk{A}(\bar{g})) \cap T = (q^\frk{A}(\bar{g}) \cap T) \cap (r^\frk{A}(\bar{g}) \cap T) \in \la G \vert_T \cup G' \ra_{\frk{A} \vert_T}
  \]

  Suppose $p = -q$.
  Then
  \[
    p^\frk{A}(\bar{g}) \cap T = -q^\frk{A}(\bar{g}) \cap T = -(q^\frk{A}(\bar{g}) \cap T) \cap T = (q^\frk{A}(\bar{g}) \cap T)' \in \la G_T \cup G' \ra_{\frk{A} \vert_T}
  \]
  where $(-)'$ is relative complementation in $\frk{A}_T$.

  Suppose $p = \exists q$.
  Then
  \begin{align*}
    p^\frk{A}(\bar{g}) \cap T &= E(q^\frk{A}(\bar{g})) \cap T \\
    &= E((q^\frk{A}(\bar{g}) \cap -T) \cup (q^\frk{A}(\bar{g}) \cap T)) \cap T\\
    &= (E(q^\frk{A}(\bar{g}) \cap -T) \cup E(q^\frk{A}(\bar{g}) \cap T)) \cap T\\
    &= E(q^\frk{A}(\bar{g}) \cap T) \cap T = E \vert_T (q^\frk{A}(\bar{g}) \cap T) \in \la G \vert_T \cup G' \ra_{\frk{A} \vert_T}
  \end{align*}
  since $-T$ is $E$-saturated.

  Suppose $p = \lozenge q$.
  Then 
  \begin{align*}
    p^\frk{A}(\bar{g}) \cap T &= R^{-1}(q^\frk{A}(\bar{g})) \cap T \\
    &= R^{-1}((q^\frk{A}(\bar{g}) \cap -T) \cup (q^\frk{A}(\bar{g}) \cap T)) \cap T\\
    &= (R^{-1}(q^\frk{A}(\bar{g}) \cap -T) \cup R^{-1}(q^\frk{A}(\bar{g}) \cap T)) \cap T\\
    &= (R^{-1}(q^\frk{A}(\bar{g}) \cap -T) \cap T) \cup (R^{-1}(q^\frk{A}(\bar{g}) \cap T) \cap T) \\
    &= U \cup R \vert_T^{-1}(q^\frk{A}(\bar{g}) \cap T) \in \la G \vert_T \cup G' \ra_{\frk{A} \vert_T}
  \end{align*}
  for some $U \in G'$, since $q^\frk{A}(\bar{g}) \cap -T \in \frk{A}'$.

  Since $\frk{A}_T \in \bb{V} \vert_T$, this establishes that $\abs{\frk{A}_T} \leq f(m + g(m))$.
  We conclude that $\abs{\frk{A}} \leq f(m + g(m)) \cdot g(m)$.
  Hence $(\bb{V} + P^0_n)_\text{SI}$ is uniformly locally finite and, by \cite[Thm.~3.7(2)]{Bez01}, $\bb{V} + P^0_n$ is locally finite.
\end{proof}

Though this result is all we will need to obtain a syntactic characterization in the next section, we provide the following alternative characterization of the variety $\bb{V} \vert_T$ -- namely, that is is the variety generated by (clopen) $Q$-clusters appearing in algebras from $\bb{V}$ (by a \textit{$Q$-cluster} we mean an equivalence class of the equivalence relation $E_Q = Q \cap Q^{-1}$).

\begin{theorem}
  Let $\bb{V} \subseteq \mpsfv$ of finite depth. Then
  \[
    \bb{V} \vert_T = \bb{Var}(\setbuilder{\frk{A} \vert_U}{\frk{A} \in \bb{V} \tand U \text{ is a clopen } Q\text{-cluster in } \frk{A}_*})
  \]
\end{theorem}

\begin{proof}
  The $(\subseteq)$ direction is immediate, since the bottom layer of any s.i. algebra from $\bb{V}$ is a $Q$-cluster (indeed, the set of $Q$-roots); hence every generator of $\bb{V} \vert_T$ is among the generators of the variety on the right.
  Conversely, take $\frk{A} \in \bb{V}$ and $U$ a clopen $Q$-cluster in the dual frame $\frk{F}$.
  $U$ must be flat; otherwise, $x R^+ y$ for some $x, y \in U$ and, by \cref{lem:depths-coincide-points}, $E(x) \, \overline{R}^+ \, E(y)$; in particular, $y \not\in Q(x)$, contrary to $U$ being a $Q$-cluster.
  Note also that this means $U \subseteq D_n(\frk{F})$ for some $n \leq d(\frk{F})$.
  Now $Q(U)$ is a closed $Q$-upset of $\frk{F}$, so $\frk{F}' = (Q(U), R \vert_U, E \vert_U)$ is the dual space of a homomorphic image $\frk{A}'$ of $\frk{A}$ (see, e.g., \cite[Thm.~3.4]{BM23}).
  We have $T(\frk{F}') = U$: certainly every element of $U$ is a $Q$-root of $\frk{F}'$; for the other inclusion, suppose $r$ is a $Q$ root of $\frk{F}'$; then $r Q x$ for some $x \in U$, so $x Q r$ and $r \in U$.
  Then $\frk{A}'$ is an s.i. algebra from $\bb{V}$ for which $\frk{A}' \vert_T = \frk{A} \vert_U$, so every generator of the variety on the right belongs to $\bb{V} \vert_T$.
\end{proof}

\begin{remark}
  \label{rem:cluster-connection}
  This previous theorem reveals a connection to the results of \cite{Shap16}.
  It follows from \cite[Cor.~4.9]{Shap16} and \cite[Thm.~5.2]{Shap16} (and is restated in \cite[Thm.~2.9]{Shap24}) that a Kripke complete logic $L$ is locally tabular iff it is of uniformly finite depth and the logic $L_C$ of clusters occurring in Kripke frames for $L$ is locally tabular.
  Indeed, if $L$ is a Kripke complete extension of $\mpsfl$, then $L$ is uniformly finite depth (height) in the sense of \cite{Shap16} iff the class of $L$-frames has bounded $Q$-depth (as $Q$ is the reachability relation in $\msf{MS4}$-frames) iff it has bounded depth (in the present sense of $R$-depth, since $d^0(\frk{F}) \leq d(\frk{F})$).
  Similarly, in this case the variety $\bb{V} \vert_T$ is generated by the algebras of \textit{all} $Q$-clusters appearing in Kripke frames for $L$ (since $\bb{V}$ is generated by the Kripke frames for $L$), and its logic coincides with $L_C$.
  Thus in this case, the criteria of \cref{thm:mpsfl-semantic-criterion} are equivalent to those of \cite{Shap16}.
\end{remark}

\subsection{A syntactic criterion for local finiteness}

\begin{definition}
  Let $\frk{A} = (B, \lozenge, \exists)$ be an $\msf{MS4}$-algebra and $\frk{F} = (X, R, E)$ its dual space.
  \begin{enumerate}
    \item Define $\diamondsuit = \lozenge \vee \exists$.
    \item Define $S = R \cup E$.
  \end{enumerate}
\end{definition}

It is straightforward to see that $\diamondsuit$ is a (possibility) operator on $B$ and $S$ is its dual relation.

The following definition is introduced (in greater generality) in \cite{Shap16} and later \cite{Shap24}.

\begin{definition}
  \label{def:rpp}
  Let $\text{RP}_m$ denote the first-order formula
  \[
    \forall x_0, \dots, x_{m+1} \left( x_0 S x_1 S \dots S x_{m+1} \to \bigvee_{i < j \leq m+1} x_i = x_j \vee \bigvee_{i < j \leq m} x_i S x_{j+1} \right)
  \]
  (the \textit{reducible path property}) and $\msf{rp}_m$ the modal formula:
  \[
    p_0 \wedge \diamondsuit(p_1 \wedge \diamondsuit(p_2 \wedge \dots \wedge \diamondsuit p_{m+1})) \to \bigvee_{i < j \leq m+1} \diamondsuit^i(p_i \wedge p_j) \vee \bigvee_{i < j \leq m} \diamondsuit^i (p_i \wedge \diamondsuit p_{j+1})
  \]
  We refer to any sequence $x_0, \dots, x_{m+1}$ of points in a frame that falsifies $\text{RP}_m$ as an \textit{irreducible path}, and we say its length is $m+1$.
  We call it \textit{proper} if there is no $j$ with $x_j (R \cap E) x_{j+1}$.
\end{definition}

\begin{remark}
  The only non-proper irreducible path is of length one, e.g. $x_0 (R \cap E) x_1$; for longer paths, such an occurrence is prohibited by irreducibility.
  This distinction will usually not be relevant, but is important for some minor technical details (e.g., \cref{lem:irreducible-path-construction}).
\end{remark}

\needspace{\baselineskip}
\begin{theorem}\
  \label{thm:rpm-frame-condition}
  \begin{enumerate}
    [label=\normalfont(\arabic*), ref = \thetheorem(\arabic*)]
    \item \label[theorem]{thm:rpm-frame-condition:1} $\msf{rp}_m$ is a canonical formula.
    \item \label[theorem]{thm:rpm-frame-condition:2} $\frk{A} \models \msf{rp}_m$ iff $\frk{A}_* \models \text{RP}_m$
  \end{enumerate}
\end{theorem}

\begin{proof}
  (1) is \cite[Prop.~3.17]{Shap24} and follows from the fact that $\msf{rp}_m$ is a Sahlqvist formula.

  (2) As stated in \cite[Prop.~3.15]{Shap24}, $\text{RP}_m$ is a first-order correspondent for $\msf{rp}_m$ on Kripke frames (and $S$ is the dual relation of $\diamondsuit$).
  But by (1), $\msf{rp}_m$ is valid on a descriptive frame iff it is valid on the underlying Kripke frame (that is, the descriptive frame condition for $\msf{rp}_m$ coincides with the Kripke frame condition; see \cite[Prop.~5.85]{BRV01}).
\end{proof}

$\msf{rp}_m$ can be thought of as a very strong form of $m$-transitivity, in that any path of length $> m$ can be ``short-circuited'' -- indeed, an irreducible path in a frame is analogous to an \textit{induced path} in graph theory, and $\msf{rp}_m$ puts a limit on the length of such paths.
Arbitrarily long irreducible paths are a ``trivial'' obstruction to local finiteness in a very general setting.
In particular, it follows from \cite[Cor.~4.9]{Shap16} and \cite[Thm.~7.3]{Shap16} that any locally finite variety of (arbitrary) Boolean algebras with operators must validate some $\msf{rp}_m$.
In particular,

\begin{theorem}
  \label{thm:rpm-necessary}
  If $\bb{V} \subseteq \bb{MS4}$ is locally finite, then $\bb{V} \models \msf{rp}_m$ for some $m < \omega$.
\end{theorem}

In \cite[Cor.~4.9,Thm.~4.10]{Bez02} it was shown that every proper subvariety of $\bb{S5}^2$ is locally finite; here we note that in fact the formulas $\msf{rp}_m$ provide a syntactic criterion for subvarieties of $\bb{S5}^2$:

\begin{theorem}
  \label{thm:s52-synactic-criterion}
  For $\bb{V} \subseteq \bb{S5}^2$, the following are equivalent:
  \begin{enumerate}
    \item $\bb{V}$ is locally finite
    \item $\bb{V} \subset \bb{S5}^2$
    \item $\bb{V} \models \msf{rp}_m$ for some $m < \omega$.
  \end{enumerate}
\end{theorem}

\begin{proof}
  The equivalence of (1) and (2) is \cite[Cor.~4.9]{Bez02}, and $(1) \Rightarrow (3)$ by \cref{thm:rpm-necessary} (\cite{Shap16}).
  Finally, for $(3) \Rightarrow (2)$, it is enough to observe that $\bb{S5}^2 \not\models \msf{rp}_m$ for any $m$.
  Indeed, $\bb{S5}^2$ contains the dual algebra of the (Kripke) frame $\frk{F} = (\omega \times \omega, E_1, E_2)$, where the equivalence classes of $E_1$ and $E_2$ are rows and columns, respectively (see \cite[Def.~3.1]{Bez02} and, e.g., \cite[Ex.~5.10]{BM23}); evidently $\frk{F}$ contains irreducible paths of arbitrary length, e.g. $(0, 0) \, S \, (1, 0) \, S \, (1, 1) \, S \, (2, 1) \dots$, and the claim follows.
\end{proof}

\cref{thm:mpsfl-semantic-criterion} shows that local finiteness of a variety of $\mpsfl$-algebras comes down to the local finiteness of a variety of $\msf{S5}^2$-algebras.
Using this, we obtain the promised syntactic criterion for $\mpsfv$:

\begin{theorem}
  \label{thm:mpsfl-syntactic-criterion}
  $\bb{V} \subseteq \mpsfv$ is locally finite iff
  \begin{enumerate}
    \item $\bb{V} \models P_n$ for some $n$ ($\bb{V}$ is finite depth)
    \item $\bb{V} \models \msf{rp}_m$ for some $m$.
  \end{enumerate}
\end{theorem}

\begin{proof}
  Suppose $\bb{V}$ is locally finite.
  Then (1) holds by the same argument as in \cref{thm:mpsfl-semantic-criterion}, and (2) holds by \cref{thm:rpm-necessary}.
  Now assume (1) and (2).
  By \cref{thm:mpsfl-semantic-criterion}, it suffices to show that $\bb{V} \vert_T$ is locally finite or, equivalently by \cref{thm:s52-synactic-criterion}, that $\bb{V} \vert_T \models \msf{rp}_m$.
  Indeed, take a s.i. $\frk{A} \in \bb{V}$ with dual space $\frk{F}$, so that $\frk{A} \vert_T$ is a generator of $\bb{V} \vert_T$.
  The dual space of $\frk{A} \vert_T$ is $\frk{F} \vert_T = (T(\frk{F}), R \vert_T, E \vert_T)$.
  Then any path $x_0 \, S \vert_T \, x_1 \, S \vert_T \dots S \vert_T \, x_{m+1}$ where $S \vert_T = R \vert_T \cup E \vert_T$ is, in particular, a path $x_0 \, S \, x_1 \, S \dots S \, x_{m+1}$ in $\frk{F}$, and hence it can be reduced (in $\frk{F})$. 
  We claim it can also be reduced in $\frk{F} \vert_T$.
  If some $x_i = x_j$ there is nothing to show; if $x_i S x_{j+1}$ for some $i < j$, then since $x_i, x_{j+1}$ both lie in the same layer $T(\frk{F})$, $x_i \, S \vert_T \, x_{j+1}$.
  This shows that any path of length $m+1$ in $\frk{F} \vert_T$ can be reduced, and hence $\frk{A} \vert_T \models \msf{rp}_m$.
  This is the case for any generator of $\bb{V} \vert_T$, so $\bb{V} \vert_T \models \msf{rp}_m$.
\end{proof}

It was pointed out by Shapirovsky and Sliusarev\footnote{via private communication} that the logic $\msf{S4}[n] \times \msf{S5}$ of products of $\msf{S4}[n]$-Kripke frames with $\msf{S5}$-Kripke frames is actually an extension of $\mpsfl$; hence \cref{thm:mpsfl-syntactic-criterion} applies to extensions of $\msf{S4}[n] \times \msf{S5}$ as well, and such an extension is locally tabular iff it contains some $\msf{rp}_m$.
In the terminology of \cite{Shap24}, this means that $\msf{S4}[n] \times \msf{S5}$ has the \textit{product rpp criterion}, as conjectured in \cite[22]{Shap24}.
In the remainder of the subsection, we will present the argument for this result as an application of \cref{thm:mpsfl-syntactic-criterion}.

\begin{definition}[{\cite[Sec.~5.1]{GKWZ03}}]
  The \textit{product} of an $\msf{S4}$-Kripke frame $(X, R)$ and an $\msf{S5}$-Kripke frame $(Y, E)$ is the frame $(X \times Y, R_1, E_2)$ where $(x, y) R_1 (x', y')$ iff $x R y \tand y = y'$, and $(x, y) E_2 (x', y')$ iff $x = x' \tand y E y'$.
  Let $\bb{S4}[n] \times \bb{S5}$ denote the variety generated by algebras of products of $\msf{S4}[n]$-Kripke frames and $\msf{S5}$-Kripke frames, which we will refer to as $(\msf{S4}[n], \msf{S5})$-product frames.
\end{definition}

The logic $\msf{S4}[n] \times \msf{S5}$ is an extension of the \textit{commutator} $[\msf{S4}, \msf{S5}]$, the logic of a commuting $\msf{S4}$ and $\msf{S5}$-operator satisfying a Church-Rosser property (\cite[Sec.~5.1]{GKWZ03}).
This commutator also corresponds precisely to the variety $\bb{MS4B}$ obtained from $\bb{MS4}$ by asserting the \textit{Barcan axiom}, as will be pointed out in \cref{def:barcan} and the following comments.
In particular, $\bb{S4}[n] \times \bb{S5} \subseteq \bb{MS4}$.
However, as was pointed out by Shapirovsky and Sliusarev, $(\msf{S4}[n], \msf{S5})$-product frames actually validate $\msf{M^{+}Cas}$, as follows from the next lemma:

\begin{lemma}
  \label{lem:kripke-frame-cas}
  Let $\frk{F} = (X, R, E)$ be an $\msf{MS4}$-Kripke frame such that $R^+$ is Noetherian (i.e., there are no proper infinite ascending $R$-chains) and all $E$-clusters are flat.
  Then $\frk{F} \models \msf{M^{+}Cas}$.
\end{lemma}

\begin{proof}
  The hypotheses allow us to mimic one direction of the proof of \cref{thm:cas-frame-condition}; indeed, suppose toward contradiction that $\frk{F} \not\models \msf{M^{+}Cas}$.
  Then take a valuation $\nu$ and $x \in X$ such that $x \not\models \msf{M^{+}Cas}$, so $x \not\models \blacksquare p$ but $x \models \blacksquare (\square (\square p \to \blacksquare p) \to \blacksquare p)$.
  Then find $x Q y$ with $y \not\models p$, and choose $z \in \qmax (R(y) \cap E(\bbrackets{-p}))$ ($R^+$ is Noetherian iff every non-empty subset contains a quasimaximal point).
  By assumption, $E(z)$ is flat.
  Since $z \in E(\bbrackets{-p})$, we may find $z E w$ with $w \not\models p$.
  Since $x Q w$ and $w \not\models \blacksquare p$, we must have $w \not\models \square(\square p \to \blacksquare p)$.
  Then we may find $w R u$ for $u \models \square p$ but $u \not\models \blacksquare p$, and further $u R v E t$ with $t \not\models p$.
  Since $z \in \qmax E(\bbrackets{-p})$ and $E(z)$ is flat, $w \in \qmax E(\bbrackets{-p})$ by \cref{lem:flat-cluster-lemmas:1}.
  Then since $w R v$ and $v \in E(\bbrackets{-p})$, $vRw$ and hence $u R w$, a contradiction since $u \models \square p$ but $w \not\models p$.
\end{proof}

\begin{lemma}
  \label{lem:product-cas}
  Let $\frk{F} = (X \times Y, R_1, E_2)$ be the product of the $\msf{S4}[n]$-Kripke frame $(X, R)$ and the $\msf{S5}$-Kripke frame $(Y, E)$.
  Then $\frk{F} \models \msf{M^{+}Cas}$.
\end{lemma}

\begin{proof}
  It is immediate from the definition that $R_1^+$ is Noetherian, since $R$ is finite depth, and that $R_1 \cap E_2$ is the diagonal relation on $X \times Y$ -- in particular, every $E_2$-cluster is flat.
  Then the claim follows from \cref{lem:kripke-frame-cas}.
\end{proof}

Certainly $\bb{S4}[n] \times \bb{S5} \models P_n$, so from \cref{thm:mpsfl-syntactic-criterion} we obtain the following:

\begin{corollary}\
  \label{cor:product-criterion}
  \begin{enumerate}
    \item $\bb{S4}[n] \times \bb{S5} \subseteq \mpsfv$
    \item $\bb{V} \subseteq \bb{S4}[n] \times \bb{S5}$ is locally finite iff $\bb{V} \models \msf{rp}_m$ for some $m$
    \item The product rpp criterion holds for $\msf{S4}[n] \times \msf{S5}$ (cf. \cite[22]{Shap24})
  \end{enumerate}
\end{corollary}

\subsection{Minimal varieties in \texorpdfstring{$\mpsfv$}{M+S4}}

In the remainder of the section, we show that $\mpsfv$ is similar to $\bb{S4}$ in that each class of varieties of depth $n$ has a unique minimal subvariety, for which we give an explicit axiomatization.
This will yield the decidability of (1) in \cref{thm:mpsfl-syntactic-criterion}.

\begin{definition}
  For $1 \leq n < \omega$, let $\frk{L}_n$ be the $\mpsfl$-frame that is an $n$-element descending chain with discrete equivalence relation, i.e. $(\set{0, \dots, n-1}, \geq, =)$.
\end{definition}

\begin{definition}
  Let
  \begin{itemize}
    \item $\msf{grz} = \square(\square(p \to \square p) \to p) \to p$
    \item $\msf{sc} = \square(\square p \to q) \vee \square(\square q \to p)$
    \item $\msf{ed} = \exists p \to \lozenge p$
  \end{itemize}
  and define
  \begin{itemize}
    \item $\bb{L}_n = \mpsfv + \msf{grz} + \msf{sc} + \msf{ed} + P_n$ (for $1 \leq n < \omega$)
    \item $\bb{L}_\omega = \mpsfv + \msf{grz} + \msf{sc} + \msf{ed}$
  \end{itemize}
\end{definition}

We recall the (descriptive) frame conditions for these formulas:

\begin{definition}[{\cite[Def.~3.5.4]{Esa19}}]
  Let $\frk{F}$ be an $\msf{MS4}$-frame and $U$ a clopen set.
  We say $x \in U$ is \textit{active} (in $U$) if there exist $y \not\in U$ and $z \in U$ such that $x R y R z$; otherwise $x$ is \textit{passive} (in $U$).
  We denote the passive points of $U$ by $\pi U$.
\end{definition}

\begin{lemma}
  \label{lem:min-variety-axiom-conditions}
  Let $\frk{F}$ be a descriptive $\msf{MS4}$-frame.
  \begin{enumerate}
    [label=\normalfont(\arabic*), ref = \thelemma(\arabic*)]
    \item \label[lemma]{lem:min-variety-axiom-conditions:1} $\frk{F} \models \msf{grz}$ iff, for each clopen $U$, $\qmax U = \max U$.
    \item \label[lemma]{lem:min-variety-axiom-conditions:2} $\frk{F} \models \msf{grz}$ iff, for each clopen $U$ and $x \in U$, there exists $y \in \pi U$ with $x R y$.
    \item \label[lemma]{lem:min-variety-axiom-conditions:3}  $\frk{F} \models \msf{sc}$ iff $R$ is \textit{strongly connected} -- that is, it satisfies the first-order condition $\forall x,y,z, x R y \wedge x R z \to y R z \vee z R y$.
    \item \label[lemma]{lem:min-variety-axiom-conditions:4}  $\frk{F} \models \msf{ed}$ iff $E \subseteq R$.
  \end{enumerate}
\end{lemma}

\begin{proof}
  (1) and (2) are characterizations of Esakia for descriptive $\msf{S4}$-frames (\cite[Thm.~3.5.5]{Esa19} and \cite[Thm.~3.5.6]{Esa19}, respectively); clearly these results hold for $\msf{MS4}$-frames as well.
  That (3) corresponds to the condition of strong connectedness on Kripke frames is well-known (\cite[Prop.~3.40]{CZ97}) and, since $\msf{sc}$ is a Sahlqvist and hence canonical formula (\cite[Def.~4.30]{BRV01} together with \cite[Prop.~5.85]{CZ97}), it corresponds to the same condition on descriptive $\msf{S4}$-frames, and hence $\msf{MS4}$-frames (see \cite[Prop.~5.85]{BRV01}).
  Finally, it is straightforward to verify that $\forall x,y, x E y \to x R y$ (that is, $E \subseteq R$) is the first-order correspondent of $\msf{ed}$ and the formula is again Sahlqvist and canonical, so this frame condition holds over descriptive frames as well, establishing (4).
\end{proof}

\begin{lemma}
  \label{lem:filtration-choice:1}
  Let $\frk{F}$ be an $\msf{MS4}$-frame with $\frk{F} \models \msf{grz}$ and $\frk{F} \models \msf{ed}$.
  For any clopen $U$ and $x \in U$, there exists $y \in \max U$ with $x R y$ and $E(y) = \set{y}$.
\end{lemma}

\begin{proof}
  Let $x \in U$.
  We have $x R y$ for $y \in \qmax U$ by \cref{thm:sees-maximal:1}.
  But $y \in \max U$ by \cref{lem:min-variety-axiom-conditions:1}.
  Now if $y E z$ for some $z$, then $y R z$ by \cref{lem:min-variety-axiom-conditions:4}, so $y=z$ and hence $E(y) = \set{y}$.
\end{proof}

\begin{lemma}
  \label{lem:filtration-choice:2}
  Let $\frk{F}$ be an $\msf{MS4}$-frame with $\frk{F} \models \msf{grz}$, and $U$ a clopen set.
  \begin{enumerate}
    \item $\max U \subseteq \pi U$
    \item For any $x \in \max U$ and $y \in X$, $x R y$ and $y R x$ implies $y = x$.
  \end{enumerate}
\end{lemma}

\begin{proof}
  (1) Let $y \in \max U$. 
  Then by \cref{lem:min-variety-axiom-conditions:2}, $y R z$ for some $z \in \pi U$, but by maximality, $z = y$; hence $y$ is passive.
  
  (2)
  Let $y \in X$ such that $x R y$ and $y R x$.
  If $y \in U$, then $y = x$ by maximality of $x$; otherwise, $y \not\in U$ which means that $y$ is an active point, contradicting (1).
\end{proof}

\begin{lemma}\
  \label{lem:minimal-variety-axiomitization}
  \begin{enumerate}
    \item For $1 \leq n < \omega$, $\bb{L}_n = \bb{Var}(\frk{L}_n^*)$
    \item $\bb{L}_\omega = \bb{Var}(\setbuilder{\frk{L}_n^*}{n < \omega})$
  \end{enumerate}
\end{lemma}

\begin{proof}
  The $\subseteq$ direction in both claims follows since each $\frk{L}_n$ clearly validates all the axioms defining $\bb{L}_n$ or $\bb{L}_\omega$.
  For the other direction of (1), it suffices to show that if $\bb{L}_n \not\models \varphi$, then $\varphi$ may be falsified on $\frk{L}_n$.
  Indeed let $\frk{F} = (X, R, E)$ be a dual frame from $\bb{L}_n$ and $\nu$ a valuation such that $\frk{F} \not\models_\nu \varphi$.
  We will use selective filtration to extract from $\frk{F}$ a frame isomorphic to $\frk{L}_m$ for $m \leq n$ that falsifies $\varphi$.
  
  Let $\Sigma$ be the set of subterms of $\varphi$.
  We will build by induction a sequence of finite subsets $X_i \subseteq X$ such that $E\vert_{X_i}$ is discrete and $R\vert_{X_i}$ is a total order (that is, a partial order in which any two points are related).
  For the base step, take $X_0 = \set{x}$ for some $x \in \max \bbrackets{-\varphi}$ with $E(x) = \set{x}$ (by \cref{lem:filtration-choice:1}, since $\bbrackets{-\varphi}$ is nonempty).
  Suppose we have built $X_i$ as described, and let $x \in X_i$.
  If $\exists \psi \in \Sigma$ and $x \models \exists \psi$, then $x \models \psi$, since $E(x) = \set{x}$; hence there is no need to add $E$-witnesses (for any points).
  Let $W^\lozenge(x) = \setbuilder{\lozenge \psi \in \Sigma}{x \models \lozenge \psi, x \not\models \psi}$
  be the $\lozenge$-subformulas in $\Sigma$ for which $x$ needs a witness.
  For each $\lozenge \psi \in W^\lozenge(x)$, we must have $x R u$ for some $u \in \bbrackets{\psi}$; so choose by \cref{lem:filtration-choice:1} a point $y(x, \lozenge\psi) \in \max \bbrackets{\psi}$ with $E(y(x, \lozenge\psi)) = \set{y(x, \lozenge\psi)}$ such that $x \, R \, y(x, \lozenge \psi)$.
  Let $W = \setbuilder{y(x, \lozenge\psi)}{\lozenge\psi \in W^\lozenge(x), x \in X_i}$ and $X_{i+1} = X_i \cup W$.
  Clearly $E\vert_{X_{i+1}}$ is discrete.
  Since $R$ is strongly connected (by \cref{lem:min-variety-axiom-conditions:2}), so is $R \vert_{X_{i+1}}$; then since $R_{i+1}$ is rooted, it follows that $R_{i+1}$ is total.
  Every point of $X_{i+1}$ is maximal in some clopen set, so by \cref{lem:filtration-choice:2}, $R\vert_{X_{i+1}}$ is antisymmetric.

  Since any two points in $X_i$ disagree on some formula in $\Sigma$ (each one is maximal in some $\bbrackets{\psi}$ for $\psi \in \Sigma$), we will eventually have $W = \varnothing$ and $X_k = X_{k+n}$ for some $k$.
  $\frk{F}_k = (X_k, R\vert_{X_k}, E\vert_{X_k})$ is a finite frame on which $R\vert_{X_k}$ is a total order and $E\vert_{X_k}$ is discrete; in particular $\frk{F}_k$ is isomorphic to some $\frk{L}_m$ and, since $d(\frk{F}) \leq n$, we must have $m \leq n$.
  Define a valuation on $\frk{F}_k$ by $\nu_k = \nu \vert_{X_k}$.
  The construction guarantees that, for any $x \in X_k$ and $\psi \in \Sigma$, $(\frk{F}, x) \models_\nu \psi$ iff $(\frk{F}_k, x) \models_{\nu_k} \psi$.
  Then $\varphi$ is falsified on $\frk{F}_k$ and hence on $\frk{L}_m$; since $\frk{L}_m$ is a generated subframe of $\frk{L}_n$, it is also falsified on $\frk{L}_n$, and we have (1).
  The same argument establishes (2) by noting in this case that $\frk{L}_m^*$ is simply among the given generating set.
\end{proof}

\begin{lemma}
  Let $\frk{A}$ be an $\mpsfl$-algebra with $d(\frk{A}) \geq n$.
  Then $\frk{L}_n$ is a subalgebra of $\frk{A}$.
\end{lemma}

\begin{proof}
  Let $\frk{A} = (B, \lozenge, \exists)$.
  By \cref{lem:depths-coincide-frames}, $d^0(\frk{A}) \geq n$ as well.
  That is, the $\msf{S4}$-algebra $\frk{A}^0 = (B^0, \lozenge^0)$ of $\exists$-fixpoints has depth $\geq n$.
  Then the Heyting algebra $O(\frk{A}^0)$ of $\square$-fixpoints also has depth $\geq n$ (that is, it satisfies exactly the same Heyting algebra identities whose G\"odel translations are the formulas $P_n$; see, e.g., \cite[Sec.~3]{Bez00}).
  By \cite[Cor.~6]{Bez00}, the Heyting algebra $H_n$ whose dual frame is the $n$-element descending chain is a Heyting subalgebra of $O(\frk{A}^0)$.
  It follows from \cite[Prop.~2.5.8,2.5.9]{Esa19} that $B(H_n)$ is a subalgebra of $B(O(\frk{A}^0))$, where $B(H)$ denotes the $\msf{S4}$-algebra obtained from the Boolean envelope of $H$ (\cite[Construction~2.5.7]{Esa19}); but since $H_n$ is finite, $B(H_n)$ is just the $\msf{S4}$-algebra with the same dual frame as $H_n$, so $B(H_n) = (\frk{L}_n^*)^0$.
  By \cite[Thm.~2.5.28]{Esa19}, $B(O(\frk{A}^0))$ is an $\msf{S4}$-subalgebra of $\frk{A}^0$.
  We have established $(\frk{L}_n^*)^0 \hookrightarrow \frk{A}^0$.
  But then since $(B^0, \lozenge^0, \exists)$ is a subalgebra of $\frk{A}$ on which $\exists$ is the identity, $\frk{L}_n^*$ itself is a subalgebra of $\frk{A}$.
\end{proof}

\begin{theorem}
  \label{thm:minimal-subvarieties}
  For $1 \leq n \leq \omega$, $\bb{L}_n$ is the unique minimal variety contained in $\setbuilder{\bb{V} \subseteq \mpsfv}{d(\bb{V}) = n}$.
\end{theorem}

\begin{proof}
  Let $\bb{V}$ be a variety with $d(\bb{V}) = n$; by \cref{lem:minimal-variety-axiomitization}, it suffices to show that $\bb{Var}(\frk{L}_n) \subseteq \bb{V}$ or, when $n = \omega$, $\bb{Var}(\setbuilder{\frk{L}_m}{m < \omega}) \subseteq \bb{V}$.
  Indeed, by definition, there is some $\frk{A} \in \bb{V}$ of depth $n$.
  By the previous lemma, $\frk{L}_n$ is a subalgebra of $\frk{A}$ or, in the case that $n = \omega$, every $\frk{L}_m$ is a subalgebra of $\frk{A}$.
  Therefore $\frk{L}_n \in \bb{V}$ or, if $n=\omega$, every $\frk{L}_m \in \bb{V}$.
\end{proof}

\begin{corollary}
  \label{cor:decidable-depth}
  Let $\bb{V} = \mpsfv + \varphi$ for some formula $\varphi$.
  It is decidable whether $\bb{V}$ is of finite depth ($\bb{V} \models P_n$ for some $n$).
\end{corollary}

\begin{proof}
  By the previous theorem, we must decide whether $\bb{V} \supseteq \bb{L}_\omega$ or, equivalently, whether $\bb{L}_\omega \models \varphi$.
  The modal logic (or equational theory) of $\bb{L}_\omega$ is decidable since it is generated by a recursively enumerable set of finite algebras. 
\end{proof}

\section{Local finiteness in \texorpdfstring{$\bb{MS4B}$}{MS4B}}
\label{sec:products}

\begin{definition}
  \label{def:barcan}
  Let $\msf{bar}$ denote the \textit{Barcan formula}
  \[
    \msf{bar} = \lozenge \exists p \to \exists \lozenge p
  \]
  and denote $\bb{MS4B} = \bb{MS4} + \msf{bar}$.
\end{definition}
 
$\msf{bar}$ is the converse of the $\msf{MS4}$ axiom $\exists \lozenge p \to \lozenge \exists p$.
The predicate analogues of these formulas correspond to a constant-domain condition on models of first-order modal logics (\cite[Sec.~3.6]{GKWZ03}).
The logic $\msf{MS4B}$ is precisely the \textit{commutator} $[\msf{S4}, \msf{S5}]$ -- the logic of a commuting $\msf{S4}$ and $\msf{S5}$-operator satisfying a Church-Rosser property (\cite[Sec.~5.1]{GKWZ03}).
By \cite[Cor.~5.10]{GKWZ03}, this coincides with the product logic $\msf{S4} \times \msf{S5}$ (although the analogous fact does \textit{not} hold for the logics $\msf{S4}[n]$ and $\msf{S5}$; see \cref{rem:not-product-matching}).
It is straightforward to see that a descriptive $\msf{MS4}$-frame $(X, R, E)$ validates $\msf{bar}$ iff $ER \subseteq RE$; hence in $\msf{MS4B}$-frames, we have $RE = ER$.

The main result of this section is to show that the reducible path property suffices to characterize local finiteness for subvarieties of $\bb{MS4B}[2]$.
This extends a result of Shapirovsky and Sliusarev in \cite[Thm.~6.10]{Shap24}, which states the product logic $\msf{S4.1}[2] \times \msf{S5}$, an extension of $\msf{MS4B}$, satisfies the \textit{product rpp criterion} (i.e., any extension is locally tabular iff it contains a \textit{product rpp formula}).
In the following section, we provide evidence that the case $n \geq 3$ remains out of reach by demonstrating a translation from subvarieties of the fusion $\bb{S5}_2$ into subvarieties of $\bb{MS4B}[3]$ that preserves and reflects local finiteness.
This translation is not novel, but rather based on an observation that the translation constructed in \cite{BM23}, landing in $\bb{MS4}[2]$, may be extended to land in $\bb{MS4B}[3]$.

\subsection{Local finiteness of \texorpdfstring{$\bb{MS4B}[2]$}{MS4B[2]}}

The proof of the main result relies essentially on the severe restrictions imposed by the Barcan axiom in depth 2, and relies on two technical lemmas.
The reader may wish to read the proof of \cref{thm:ms4b-syntactic-criterion} first to see the high-level argument.

\begin{lemma}
  \label{lem:cluster-commute-condition}
  Let $\frk{F} = (X, R, E)$ be a $\msf{MS4}$-frame.
  Let $x_0, y_0 \in X$ with $x_0 R y_0$.
  \begin{enumerate}
    [label=\normalfont(\arabic*), ref = \thelemma(\arabic*)]
    \item \label[lemma]{lem:cluster-commute-condition:1} For all $x \in E(x_0)$, there exists $y \in E(y_0)$ with $xRy$.
    \item \label[lemma]{lem:cluster-commute-condition:2} $\frk{F} \models \msf{bar}$ iff: for all $y \in E(y_0)$, there exists $x \in E(x_0)$ with $x R y$.
  \end{enumerate}
  and hence both properties hold in $\msf{MS4B}$-frames.
\end{lemma}

\begin{proof}
  Both assertions are trivial if $x_0 E y_0$, so suppose not.
  For (1), suppose for some $x$ there was no such $y$.
  Then $x E x_0 R y_0$, but there is no $y$ such that $x R y E y_0$, contradicting the frame condition for the $\msf{MS4}$ axiom $RE \subseteq ER$.

  For (2), suppose for some $y$ there was no such $x$.
  Then $x_0 R y_0 E y$, but there is no $x$ with $x_0 E x R y$, contradicting the frame condition for the Barcan axiom $ER \subseteq RE$.
\end{proof}

\begin{definition}
  Let $\frk{F} = (X, R, E)$ be an $\msf{MS4}$-frame.
  We refer to equivalence classes of $E$ as \textit{$E$-clusters}.
  Define $E_R = R \cap R^{-1}$ and $E_Q = Q \cap Q^{-1}$.
  We refer to equivalence classes of $E_R$ (resp. $E_Q$) as \textit{$R$-clusters} (resp. \textit{$Q$-clusters}).
\end{definition}

Recall (\cref{thm:si-simple-ms4:2}) that an $\msf{MS4}$-algebra $\frk{A}$ is simple iff its dual space is a $Q$-cluster.

\begin{lemma}
  \label{lem:irreducible-path-construction}
  Suppose $\frk{F}$ is a simple $\msf{MS4B}[2]$-frame and, for some $m \geq 1$, both of the following hold:
  \begin{enumerate}
    \item $\abs{X/E} > m$
    \item $\abs{X/E_R} > 2m$
  \end{enumerate}
  Then $\frk{F}$ contains a proper irreducible path of length strictly greater than $m$.
\end{lemma}

\begin{proof}
  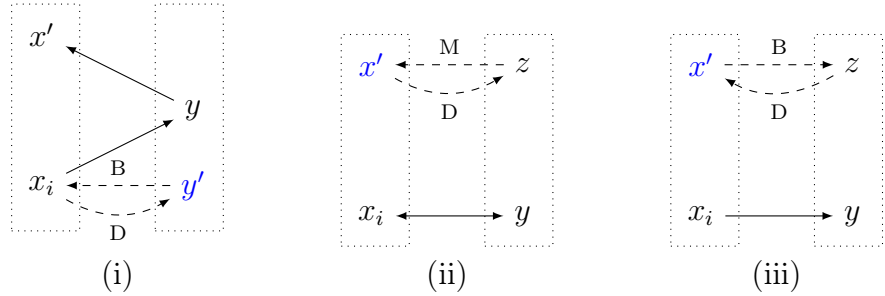
\begin{figure}[h]
    \centering
    \begin{tikzpicture}[
        scale=1,
      dot/.style={circle,fill=black,minimum size=6pt,inner sep=0}
    ]
    \node (xi) at (0, 0) {$x_i$};
    \node (y) at (2, 1) {$y$};
    \node (x') at (0, 2) {$x'$};
    \node[blue] (y') at (2, 0) {$y'$};

    \draw[-latex] (xi) -- (y);
    \draw[-latex] (y) -- (x');
    \draw[dashed,-latex] (y') -- node[midway,above] {\tiny B} (xi);
    \draw[dashed,-latex,out=330,in=210] (xi) to node[midway,below] {\tiny D} (y');

    \draw[dotted] (-0.4, -0.6) rectangle (0.5, 2.4);
    \draw[dotted] (1.5, -0.6) rectangle (2.4, 2.4);

    \node at (1, -1.2) {(i)};
    \end{tikzpicture}
    \hspace{0.5in}
    \begin{tikzpicture}[
      scale=1,
      dot/.style={circle,fill=black,minimum size=6pt,inner sep=0}
    ]
    \node (xi) at (0, 0) {$x_i$};
    \node (y) at (2, 0) {$y$};
    \node[blue] (x') at (0, 2) {$x'$};
    \node (z) at (2, 2) {$z$};

    \draw[latex-latex] (xi) -- (y);
    \draw[dashed,-latex] (z) -- node[midway,above] {\tiny M} (x');
    \draw[dashed,-latex,out=330,in=210] (x') to node[midway,below] {\tiny D} (z);

    \draw[dotted] (-0.4, -0.4) rectangle (0.5, 2.4);
    \draw[dotted] (1.5, -0.4) rectangle (2.4, 2.4);

    \node at (1, -0.8) {(ii)};
    \end{tikzpicture}
    \hspace{0.5in}
    \begin{tikzpicture}[
      scale=1,
      dot/.style={circle,fill=black,minimum size=6pt,inner sep=0}
    ]
    \node (xi) at (0, 0) {$x_i$};
    \node (y) at (2, 0) {$y$};
    \node[blue] (x') at (0, 2) {$x'$};
    \node (z) at (2, 2) {$z$};

    \draw[-latex] (xi) -- (y);
    \draw[dashed,-latex] (x') -- node[midway,above] {\tiny B} (z);
    \draw[dashed,-latex,out=210,in=330] (z) to node[midway,below] {\tiny D} (x');

    \draw[dotted] (-0.4, -0.4) rectangle (0.5, 2.4);
    \draw[dotted] (1.5, -0.4) rectangle (2.4, 2.4);

    \node at (1, -0.8) {(iii)};
    \end{tikzpicture}
    \caption{
      The key situations in \cref{lem:irreducible-path-construction}.
      The dashed arrows are obtained by an appeal to the Barcan axiom (B), the $\msf{MS4}$-axiom (M), or the depth of $\frk{F}$ (D).
      The choice of $x_{i+1}$ is depicted in blue.
    }
    \label{fig:path-construction-key-situations}
  \end{figure}

  Since every $E_R$-cluster is contained in either $D_1(\frk{F})$ or $D_2(\frk{F})$, it must be the case that one of $D_1(\frk{F})$ or $D_2(\frk{F})$ contains strictly more than $m$ $E_R$-clusters.
  Call that layer $L$, and note that $\abs{L/E_R} > m$.
  Choose any point $x_0 \in L$.
  We will build a \textit{proper} irreducible path (see \cref{def:rpp}) by induction, starting at $x_0$ and entirely contained within $L$.
  Suppose we have constructed $x_0 \, S \, x_1 \, S \dots S \, x_i$ a proper irreducible path for $i \leq m$ with all points contained in $L$.
  
  Suppose that $i=0$ or $x_{i-1} E x_i$.
  There are at most $i$ distinct $E$-clusters represented among $x_0, \dots, x_i$ (since certainly at least two are in the same cluster, at any instance where $x_j E x_{j+1}$).
  Since $i \leq m$ and $\abs{X/E} > m$, there must be a point $z$ so that $\neg(x_j E z)$ for any $j \leq i$.
  Since $x_i$ is a $Q$-root of $\frk{F}$ (as is every point, since $\frk{F}$ is a $Q$-cluster), $x_i Q z$; that is, $x_i \, R \, y$ for some $y \in E(z)$.
  If $L = D_1(\frk{F})$, then necessarily $x_i \, E_R \, y$ and $y \in L$; then we may take $x_{i+1} = y$.
  If $L = D_2(\frk{F})$ then, since $y$ is also a $Q$-root, $y R x'$ for $x' \in E(x_i)$.
  Then by \cref{lem:cluster-commute-condition:2}, there is some $y' \in E(y)$ such that $y' R x_i$; since $x_i \in D_2(\frk{F})$ and $\frk{F}$ is of depth 2, this is only possible if $y' \in D_2(\frk{F})$ (see \cref{fig:path-construction-key-situations} (i)).
  Hence $x_i \, E_R \, y'$, and we may take $x_{i+1} = y'$.
  In either case, $x_{i+1}$ was chosen so that $\neg(x_j E x_{i+1})$ for any $j \leq i$.
  It also cannot be the case that $x_j R x_{i+1}$ for any $j < i$, since then $x_j \, E_R \, x_i$; in this case, either $j \leq i-2$ and the original path was reducible, or $j = i - 1$, $x_{i-1} (R \cap E) x_i$, and the original path was not proper.
  Then $x_0 \, S \dots S \, x_{i+1}$ is a proper irreducible path contained in $L$.
  
  Suppose $x_{i-1} R x_i$.
  Similarly there are at most $i$ $E_R$-clusters represented among $x_0,\dots,x_i$.
  Since $i \leq m$ and $\abs{L/E_R} > m$, there is a point $z \in L$ such that $\neg(x_j \, E_R \, z)$ for any $j \leq i$.
  Since $x_i$ is a $Q$-root, $x_i R y$ for some $y \in E(z)$.
  Suppose $L = D_1(\frk{F})$.
  Then $x_i \, E_R \, y$ and, by \cref{lem:cluster-commute-condition:1}, there is some $x' \in E(x_i)$ with $z R x'$ and, since both are in $D_1(\frk{F})$, $x' \, E_R \, z$ (see \cref{fig:path-construction-key-situations}, (ii)).
  Then we may take $x_{i+1} = x'$.
  Suppose $L = D_2(\frk{F})$.
  Then by \cref{lem:cluster-commute-condition:2}, there is $x' \in E(x_i)$ so that $x' R z$ and, since both are in $D_2(\frk{F})$ and $\frk{F}$ is depth 2, $x' \, E_R \, z$ (see \cref{fig:path-construction-key-situations}, (iii)).
  Then we may take $x_{i+1} = x'$ again.
  In either case, $x_{i+1}$ was chosen so that $\neg(x_j \, E_R \, x_{i+1})$ for any $j \leq i$.
  It also cannot be the case that $x_j E x_{i+1}$ for any $j < i$ since then $x_j E x_{i}$, which is a contradiction by the same argument as the previous case.
  Then $x_0 S \dots S x_{i+1}$ is a proper irreducible path contained in $L$.
\end{proof}

The next lemma shows that s.i. $\msf{MS4B}[2]$-frames that are \textit{not} simple are actually $\mpsfl$-frames.

\begin{lemma}
  \label{lem:ms4b2-non-simple-cas}
  Suppose $\frk{F}$ is a s.i. $\msf{MS4B}[2]$-frame that is \textit{not} simple.
  Then $\frk{F} \models \msf{M^{+}Cas}$.
\end{lemma}

\begin{proof}
  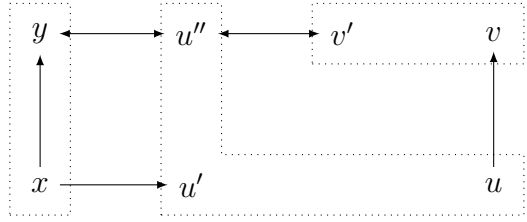
\begin{figure}[h]
    \centering
    \begin{tikzpicture}[
        scale=1,
      dot/.style={circle,fill=black,minimum size=6pt,inner sep=0}
    ]
    \node (x) at (0, 0) {$x$};
    \node (y) at (0, 2) {$y$};
    \node (u') at (2, 0) {$u'$};
    \node (u) at (6, 0) {$u$};
    \node (u'') at (2, 2) {$u''$};
    \node (v') at (4, 2) {$v'$};
    \node (v) at (6, 2) {$v$};

    \draw[-latex] (x) -- (y);
    \draw[-latex] (x) -- (u');
    \draw[latex-latex] (y) -- (u'');
    \draw[latex-latex] (u'') -- (v');
    \draw[-latex] (u) -- (v);

    \draw[dotted] (-0.4, -0.4) rectangle (0.4, 2.4);
    \draw[dotted] (1.6, -0.4) -- (1.6, 2.4) -- (2.4, 2.4) -- (2.4, 0.4) -- (6.4, 0.4) -- (6.4, -0.4) -- (1.6, -0.4);
    \draw[dotted] (3.6, 1.6) rectangle (6.4, 2.4);
    \end{tikzpicture}
    \caption{
      The situation in the proof of \cref{lem:ms4b2-non-simple-cas}.
    }
    \label{fig:ms4b2-cas-key-situation}
  \end{figure}

  $\frk{F}$ has a non-empty set of $Q$-roots $T(\frk{F})$.
  Since $\frk{F}$ is not simple, we have $d^0(\frk{F}) = 2$.
  In this situation, we will show that $T(\frk{F}) = D_2(\frk{F})$.
  For $T(\frk{F}) \subseteq D_2(\frk{F})$, suppose $x \in D_1(\frk{F})$ is a $Q$-root of $\frk{F}$; we will argue that every point of $\frk{F}$ is a $Q$-root, contradicting $d^0(\frk{F}) = 2$.
  Indeed, let $y$ be any other point in $\frk{F}$.
  Then $x R y' E y$ for some $y'$ and, since $x \in D_1(\frk{F})$, $x \, E_R \, y'$.
  Then evidently $y Q x$ and $y$ is also a $Q$-root.
  To show $D_2(\frk{F}) \subseteq T(\frk{F})$, take any $y \in D_2(\frk{F})$.
  There is a $Q$-root $x$ necessarily in $D_2(\frk{F})$ and hence $x R y' E y$ for some $y'$.
  By \cref{lem:cluster-commute-condition:2}, there is $x' \in E(x)$ such that $x' R y$ and, since $y \in D_2(\frk{F})$, so is $x'$ and $x' \, E_R \, y$.
  Then again $y Q x$ and $y$ is also a $Q$-root.

  Now assume that $\frk{F} \not\models \msf{M^{+}Cas}$.
  By \cref{thm:cas-frame-condition}, there is a non-flat cluster in $\frk{F}$.
  That is, there are $x, y$ so that $x R^+ y$ and $E(x) = E(y)$ (the following argument is depicted in \cref{fig:ms4b2-cas-key-situation}).
  It must be the case that $x \in D_2(\frk{F})$ and $y \in D_1(\frk{F})$, so $x$ is a $Q$-root.
  Now take any $u, v$ with $u R v$, so that $E(u) \, \overline{R} \, E(v)$.
  We will show that $E(v) \, \overline{R} \, E(u)$, contradicting $d^0(\frk{F}) = 2$.
  Now $x R u'$ for some $u' \in E(u)$.
  By \cref{lem:cluster-commute-condition:1}, $y R u''$ for some $u'' \in E(u)$ and, since $y \in D_1(\frk{F})$, $y \, E_R \, u''$.
  By \cref{lem:cluster-commute-condition:1} again, since $u R v$, we must have $u'' R v'$ for $v' \in E(v)$.
  But then since $u'' \in D_1(\frk{F})$, $u'' \, E_R \, v'$, and we conclude $E(v) \, \overline{R} \, E(u)$.
\end{proof}

We are finally ready to prove the main result:

\begin{theorem}
  \label{thm:ms4b-syntactic-criterion}
  $\bb{V} \subseteq \bb{MS4B}[2]$ is locally finite iff $\bb{V} \models \msf{rp}_m$ for some $m < \omega$.
\end{theorem}

\begin{proof}
  By \cref{thm:rpm-necessary}, it suffices to show $(\Leftarrow)$, so suppose $\bb{V} \models \msf{rp}_m$ for some fixed $m$.
  Let $\frk{A} \in \bb{V}$ be a s.i. algebra generated by $G = \set{g_1, \dots, g_n}$, and let $\frk{F} = (X, R, E)$ be its dual space.
  
  Suppose first that $\frk{A}$ is simple ($d^0(\frk{A}) = 1$).
  By \cref{lem:irreducible-path-construction}, at least one of the following must hold in $\frk{F}$ (if both are falsified, there is an irreducible path of length strictly greater than $m$, contradicting the frame condition for $\msf{rp}_m$):
  \begin{enumerate}
    \item $\abs{X/E} \leq m$
    \item $\abs{X/E_R} \leq 2m$
  \end{enumerate}
  Assuming (1), we have that $\abs{B^0} \leq 2^m$, where we recall $B^0$ is the set of $\exists$-fixpoints of $\frk{A}$; hence $\frk{A}$ is generated as an $\msf{S4}[2]$-algebra by $G \cup B^0$, and $\abs{\frk{A}} \leq f_1(n + 2^m)$, where $f_1 : \omega \to \omega$ is a function bounding the size of an $n$-generated algebra from the locally finite variety $\bb{S4}[2]$.
  Assuming (2), we have $\abs{K} \leq 2^{2m}$, where $K$ is the set of $\lozenge$-fixpoints of $\frk{A}$; this follows since every $\lozenge$-fixpoint of $\frk{A}$ is a clopen $R$-downset of $\frk{F}$, and every such $R$-downset is a union of $E_R$-clusters.
  In this case, $\frk{A}$ is generated as an $\msf{S5}$-algebra by $G \cup K$, and $\abs{\frk{A}} \leq f_2(n + 2^{2m})$, where $f_2 : \omega \to \omega$ is a function bounding the size of an $n$-generated algebra from the locally finite variety $\bb{S5}$.

  Now suppose $\frk{A}$ is not simple ($d^0(\frk{A}) = 2$).
  By \cref{lem:ms4b2-non-simple-cas}, we have $\frk{A} \in \mpsfv[2] + \msf{rp}_m$.
  Then $\abs{\frk{A}} \leq f_3(n)$, where $f_3 : \omega \to \omega$ is a function bounding the size of an $n$-generated algebra from $\mpsfv[2] + \msf{rp}_m$, which is locally finite by \cref{thm:mpsfl-syntactic-criterion}.

  We conclude that 
  \[
    \abs{\frk{A}} \leq \max(f_1(n + 2^m), f_2(n + 2^{2m}), f_3(n))
  \]
  hence $\bb{V}_{\text{SI}}$ is uniformly locally finite, and $\bb{V}$ is locally finite, by \cite[Thm.~3.7(2)]{Bez01}.
\end{proof}

\subsection{A translation of \texorpdfstring{$\bb{S5}_2$}{S52} into \texorpdfstring{$\bb{MS4B}[3]$}{MS4B[3]}}
\label{sec:translation}

In this section we will demonstrate a translation of subvarieties of $\bb{S5}_2$ into subvarieties of $\bb{MS4B}[3]$ that preserves and reflects local finiteness.
We briefly recall the relevant definitions, which appear in more detail in \cite{BM23}.

\begin{definition}[{\cite[Sec.~3.1]{GKWZ03}}]
  The logic $\msf{S5}_2$ is the {\em fusion} $\msf{S5} \oplus \msf{S5}$, i.e.~the smallest normal modal
  logic in the language with two modalities $\exists_1$ and $\exists_2$ containing the $\msf{S5}$ axioms for each $\exists_i$ (and no other axioms). 
\end{definition}

Algebraic models of $\msf{S5}_2$ are triples $(B, \exists_1, \exists_2)$, where $(B, \exists_i)$ is an $\msf{S5}$-algebra for $i = 1,2$; we call such algebras {\em $\msf{S5}_2$-algebras}. 
J\'{o}nsson--Tarski duality specializes to yield that the corresponding variety $\bb{S5}_2$ is dually equivalent to the following category of descriptive frames:

\begin{definition} 
  A {\em descriptive $\msf{S5}_2$-frame} is a triple $\frk{F} = (X, E_1, E_2)$ where $X$ is a Stone space and $E_1,E_2$ are two continuous equivalence relations on $X$.
\end{definition}

In \cite[Sec.~6]{BM23}, we gave a translation of subvarieties of $\bb{S5}_2$ into subvarieties of $\bb{MS4_S}[2]$ that preserves and reflects local finiteness.
Here, we will demonstrate that the same construction yields a translation into subvarieties of $\bb{MS4B_S}[3]$ (here $\bb{MS4B_S} = \bb{MS4_S} + \msf{bar}$, see \cref{thm:ms4s-semisimple}).
This demonstrates that the method of proof used in \cref{thm:ms4b-syntactic-criterion} cannot be extended to depth $\geq 3$.
The essential ingredient of the translation is the following construction on the level of individual frames:

\begin{figure}[t]
  \centering
  \begin{tikzpicture}[
    scale=1.5,
    dot/.style={circle,fill=black,minimum size=6pt,inner sep=0}
  ]
  \node (A) [dot] at (0, 1) {};
  \node (B) [dot] at (1, 1) {};
  \node (C) [dot] at (1, 0) {};
  \node (D) [dot] at (2, 0) {};
  \node (E) [dot] at (2, 1) {};
  \node (F) [dot] at (3, 1) {};
  \node (G) [dot] at (3, 0) {};
  \node (H) [dot] at (4, 0) {};

  \draw (A) -- (B);
  \draw[blue,thick] (B) -- (C);
  \draw (C) -- (D);
  \draw[blue,thick] (D) -- (E);
  \draw (E) -- (F);
  \draw[blue,thick] (F) -- (G);
  \draw (G) -- (H);

  \node at (4.5, 0.5) {$\dots$};
  \node [dot] at (5, 0.5) {};

  \node  at (-1, 0.5) {$\frk{F}$};
  \end{tikzpicture}
  
  \bigskip

  \begin{tikzpicture}[
    scale=1.25,
    dot/.style={circle,fill=black,minimum size=6pt,inner sep=0}
  ]

  \foreach \x in {0,...,4,6} {
    \draw[blue,thick,fill=blue!20] (\x - 0.5, -1.75) -- (\x - 0.5, 1.85) -- (\x + 0.25, 2.65) -- (\x + 0.25, -1) --cycle;
  }

  \node (A2) [dot] at (0, 1) {};
  \node (B2) [dot] at (1, 1) {};
  \node (C2) [dot] at (0.75, 0) {};
  \node (D2) [dot] at (1.75, 0) {};
  \node (E2) [dot] at (2, 1) {};
  \node (F2) [dot] at (3, 1) {};
  \node (G2) [dot] at (2.75, 0) {};
  \node (H2) [dot] at (3.75, 0) {};
  \node (I2) [dot] at (5.75, 0.5) {};

  \node (A1) [dot] at (0, 2) {};
  \node (B1) [dot] at (1, 2) {};
  \node (C1) [dot] at (2, 2) {};
  \node (D1) [dot] at (3, 2) {};
  \node (E1) [dot] at (4, 2) {};
  \node (F1) [dot] at (6, 2) {};

  \draw (A1) -- (E1) -- (4.5, 2);
  \draw (5.5, 2) -- (F1);
  \node at (5, 2) {$\dots$};

  \node (A3) [dot] at (0, -1) {};
  \node (B3) [dot] at (1, -1) {};
  \node (C3) [dot] at (2, -1) {};
  \node (D3) [dot] at (3, -1) {};
  \node (E3) [dot] at (4, -1) {};
  \node (F3) [dot] at (6, -1) {};

  \draw (A3) -- (E3) -- (4.5, -1);
  \draw (5.5, -1) -- (F3);
  \node at (5, -1) {$\dots$};

  \draw (A2) -- (B2);
  \draw (C2) -- (D2);
  \draw (E2) -- (F2);
  \draw (G2) -- (H2);

  \draw[-latex] (A2) -- (A1);
  \draw[-latex] (B2) -- (B1);
  \draw[-latex] (E2) -- (C1);
  \draw[-latex] (F2) -- (D1);
  \draw[-latex] (C2) to[out=90, in=250] (B1);
  \draw[-latex] (D2) to[out=90, in=250] (C1);
  \draw[-latex] (G2) to[out=90, in=250] (D1);
  \draw[-latex] (H2) to[out=90, in=250] (E1);
  \draw[-latex] (I2) to[out=90, in=250] (F1);

  \draw[-latex] (A3) -- (A2);
  \draw[-latex] (B3) -- (B2);
  \draw[-latex] (C3) -- (E2);
  \draw[-latex] (D3) -- (F2);
  \draw[-latex] (B3) to[out=120, in=270] (C2);
  \draw[-latex] (C3) to[out=120, in=270] (D2);
  \draw[-latex] (D3) to[out=120, in=270] (G2);
  \draw[-latex] (E3) to[out=120, in=270] (H2);
  \draw[-latex] (F3) to[out=120, in=270] (I2);

  \node at (5, 1) {$\dots$};

  \node at (-1.5, 0.5) {$T(\frk{F})$};
  \end{tikzpicture}
  \caption{
    Constructing a $\msf{MS4B_S}[3]$-frame from an $\msf{S5}_2$-frame.
    In the $\msf{S5}_2$-frame on top, $E_1$-clusters are horizontal lines while $E_2$-clusters are blue vertical lines.
    In the $\msf{MS4B_S}[3]$-frame on the right, $R$-clusters are horizontal lines, proper $R$-arrows are drawn with arrowheads, and $E$-clusters are given by the blue rectangles.
  }
  \label{fig:translation}
\end{figure}
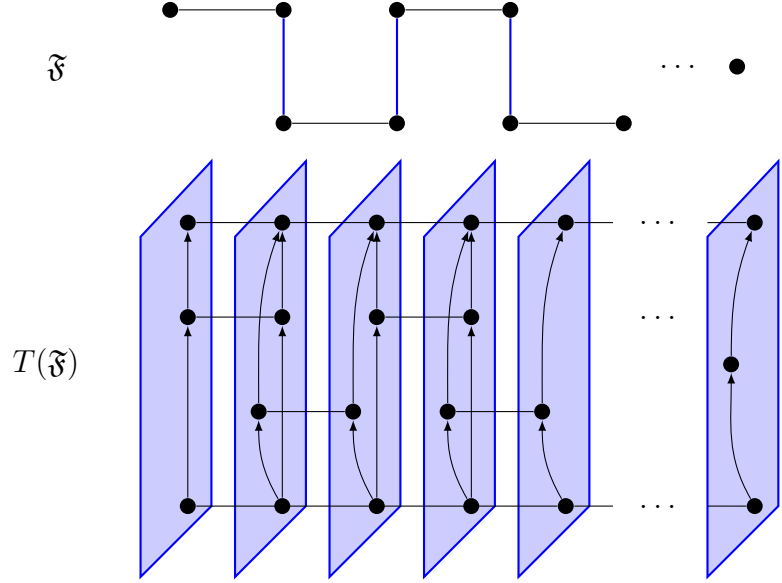

\begin{construction}
  \label{con:translation-frame}
  Let $\frk{F} = (X, E_1, E_2)$ be a descriptive $\msf{S5}_2$-frame. We let $\mca{L} = X / E_2$ be the quotient space and $\pi_\mca{L} : X \to \mca{L}$ the quotient map.
  Let $\mca{L}^T = \setbuilder{\alpha^T}{\alpha \in \mca{L}}$ and $\mca{L}_B = \setbuilder{\alpha_B}{\alpha \in \mca{L}}$ denote the two copies of $\mca{L}$ in the disjoint union $\mca{L} \sqcup \mca{L}$.
  We set $T(\frk{F}) = (Y, R, E)$, where
  \begin{itemize}
      \item $Y = X \cup \mca{L}^T \cup \mca{L}_B$ is the (disjoint) union of $X$ and the two copies of $\mca{L}$,
      \item $E$ is the smallest equivalence on $Y$ containing $E_2$ along with
      \[
      \setbuilder{(x, \alpha^T)}{x \in X, \alpha \in \mca{L}, x \in \alpha} \cup \setbuilder{(x, \alpha_B)}{x \in X, \alpha \in \mca{L}, x \in \alpha}
      \]
      \item $R = E_1 \cup (\mca{L}_B \times Y) \cup (X \times \mca{L}^T) \cup (\mca{L}^T \times \mca{L}^T)$.
  \end{itemize}
\end{construction}

This construction is depicted in \cref{fig:translation}.
The definition of $R$ makes $\mca{L}^T$ a cluster that is seen by every point (a  ``top rail''), and $\mca{L}_B$ a cluster that sees every point (a ``bottom rail'').

\begin{lemma}
  $T(\frk{F})$ is a descriptive $\msf{MS4B_S}[3]$-frame.
\end{lemma}

\begin{proof}
  The details concerning the topology and continuity of $R$ and $E$ are the same as in \cite[Lem.~6.10]{BM23}.
  Here we focus on the new additions.

  Indeed, observe that $Q = ER$ is the total relation on $T(\frk{F})$: For each $x \in X$ we have $\mca{L}^T \subseteq R(x)$, and for every other point $y$, we have $y \, E \, E_2(y)^T \in \mca{L}^T$.
  On the other hand, $RE$ is also total: For each $x \in X$, we have $x\, E \, E_2(x)_B$, and $R(E_2(x)_B) = Y$.
  Then we have $ER = RE = Y^2$ and, in particular, $T(\frk{F})$ is a $Q$-cluster.
  Evidently $T(\frk{F})$ is of depth 3.
  Hence $T(\frk{F})$ is an $\msf{MS4B_S}[3]$-frame.
\end{proof}

We extend our translation on frames to one on varieties of algebras.
For $\frk{A} \in \bb{S5}_2$, let $T(\frk{A}) = (T(\frk{A}_*))^*$.
For a variety $\bb{V} \subseteq \bb{S5}_2$, let $\bb{V}_\text{SI}$ be the class of s.i.~$\bb{V}$-algebras, and define $T(\bb{V})$ to be the variety generated by the translations of the s.i.~members:
\[
T(\bb{V}) = \bb{Var}(\setbuilder{T(\frk{A})}{\frk{A} \in \bb{V}_\text{SI}}).
\]
It follows from \cref{con:translation-frame} that $T(\bb{V})$ is generated by a class of $\msf{MS4B_S}[3]$-algebras, so $T(\bb{V}) \subseteq \bb{MS4B_S}[3]$.

The proofs in \cite[Sec.~6]{BM23} may be adapted with only superficial  alterations to yield:

\begin{theorem}
  $\bb{V} \subseteq \bb{S5}_2$ is locally finite iff $T(\bb{V}) \subseteq \bb{MS4B_S}[3]$ is locally finite.
\end{theorem}

\begin{remark}
  \label{rem:not-product-matching}
  Whereas $[\msf{S4}, \msf{S5}]$ coincides with the product logic $\msf{S4} \times \msf{S5}$, it is \textit{not} the case that $[\msf{S4}[n], \msf{S5}] = \msf{MS4B}[n]$ coincides with $\msf{S4}[n] \times \msf{S5}$:
  This can be seen by comparing \cref{cor:product-criterion} with \cref{con:translation-frame}, which provides an example of an $\msf{MS4B_S}[3]$-frame that is not an $\mpsfl$-frame. 
  This fact is also observed and pointed out in \cite[Rem.~6.4]{Shap24}.
\end{remark}

\section{Conclusion}
\label{sec:conclusion}

This paper continues the study of local finiteness in $\bb{MS4}$ undertaken in \cite{BM23}, which contains mostly negative results demonstrating that many properties of $\bb{S4}$ fail in the monadic setting.
On the other hand, this paper shows that $\mpsfv$ is much better behaved and shares many similarities with $\bb{S4}$.
Indeed, \cref{thm:mpsfl-syntactic-criterion} can be viewed as an analogue of the Segerberg--Maksimova theorem, stating that $\bb{V} \subseteq \mpsfv$ is locally finite iff it is of ``finite depth'' and ``finite width'', while \cref{thm:minimal-subvarieties} shows that there is a natural minimal subvariety contained in each class of subvarieties of a given depth.
Along with \cite{BBI23}, which shows that the Casari axiom is necessary to obtain a provability interpretation of monadic intuitionistic propositional logic, this contributes to evidence $\mpsfl$ is an interesting and well-motivated extension of $\msf{MS4}$ to study.

We end with some questions and directions for future work.
In light of \cref{cor:decidable-depth}, an obvious question is whether condition (2) in \cref{thm:mpsfl-syntactic-criterion} is also decidable, which would imply that local finiteness is decidable for subvarieties of $\mpsfv$.
The reducible path property $\msf{rp}_m$ is itself of significant interest.
\cite[Cor.~4.9]{Shap16} and \cite[Thm.~7.3]{Shap16} demonstrates that this is a necessary condition for local finiteness in any modal logic.
Indeed, the failure of $\msf{rp}_m$ is an ``trivial'' obstruction to local finiteness, and it would be interesting to investigate other families of logics for which it is sufficient.
This issue seems quite subtle, as an example of Makinson (\cite{Mak81}) shows that $\msf{rp}_m$ is not sufficient even for the (unimodal) logic $\msf{KTB}$.
Finally, \cite[Q.~5.13]{BM23} asks if $D_1(\frk{F})$ or more generally any $D_n(\frk{F})$ is clopen in a finitely generated $\msf{MS4}$-frame.
Since this holds for $\msf{S4}$-frames, there is some reason to expect that it might hold for $\mpsfl$-frames, but the question remains unresolved even in this setting; the only result we have obtained in this direction is \cref{lem:si-bottom-clopen}, which shows that the bottom layer of an s.i. $\mpsfl$-frame is clopen.

\section{Acknolwedgements}

I am thankful to Ilya Shapirovsky, Vladislav Sliusarev, and Guram Bezhanishvili for their observations and a valuable discussion resulting in the inclusion of \cref{rem:cluster-connection}, \cref{lem:kripke-frame-cas}, \cref{lem:product-cas}, and \cref{cor:product-criterion}.

\printbibliography

\end{document}